\documentclass[11pt]{article}
\usepackage{ytableau}
\usepackage{authblk}
\usepackage{amssymb}
\usepackage{amsmath, epsfig, cite,colordvi,xcolor}
\usepackage{booktabs}
\usepackage[margin=1in]{geometry}
\usepackage{mathrsfs}
\usepackage{stmaryrd}
\usepackage{enumerate}
\usepackage{epstopdf}
\usepackage{amsfonts,amsmath,amssymb,amscd,bm}
\usepackage{shadow}
\usepackage{graphicx}
\usepackage{color}
\usepackage{longtable}
\usepackage{pstricks,multido}
\usepackage[colorlinks,urlcolor=purple,linkcolor=red,anchorcolor=green,citecolor=blue]{hyperref}
\usepackage{qtree}
\allowdisplaybreaks
\newtheorem{thm}{Theorem}[section]
\newtheorem{lem}[thm]{Lemma}

\newtheorem{cor}[thm]{Corollary}
\newtheorem{remark}{Remark}[section]
\newtheorem{claim}{Claim}[section]
\newtheorem{example}{Example}[section]

\numberwithin{figure}{section}
\numberwithin{table}{section}

\def\qed{\hfill \rule{4pt}{7pt}}

\def\pf{\noindent {\it{Proof.} \hskip 2pt}}

\numberwithin{equation}{section}

\def\f{\mathcal{F}}
\begin{document}

\begin{center}
{\Large\bf Combinatorial proofs for identities related to
generalizations\\[6pt] of the mock theta functions $\omega(q)$ and $\nu(q)$ }
\end{center}

\begin{center}

Frank Z.K. Li$^{1}$
and Jane Y.X. Yang$^{2}$

$^{1}$Center for Combinatorics, LPMC \\
Nankai University, Tianjin 300071, P.R. China\\[6pt]

$^{2}$School of Science \\
Chongqing University of Posts and Telecommunications\\
Chongqing 400065, P.R. China

emails:
zhkli@mail.nankai.edu.cn, yangyingxue@cqupt.edu.cn.

\end{center}

\begin{abstract}
The two partition functions $p_\omega(n)$ and $p_\nu(n)$ were introduced by Andrews, Dixit and Yee,
which  are related to the third order mock theta functions $\omega(q)$ and $\nu(q)$, respectively.
Recently, Andrews and Yee analytically studied two  identities
that connect the refinements of $p_\omega(n)$ and $p_\nu(n)$ with the  generalized bivariate mock theta functions $\omega(z;q)$ and $\nu(z;q)$, respectively.
However, they stated
these identities cried out for bijective proofs.
In this paper, we  first define the  generalized trivariate   mock theta functions $\omega(y,z;q)$ and $\nu(y,z;q)$.
Then by utilizing  odd Ferrers graph,
we obtain  certain identities concerning to $\omega(y,z;q)$ and $\nu(y,z;q)$,
which extend some early results  of Andrews that are  related to $\omega(z;q)$ and $\nu(z;q)$.
In virtue of the combinatorial interpretations that arise from the identities involving $\omega(y,z;q)$ and $\nu(y,z;q)$,
we finally present bijective proofs  for the two identities of Andrews-Yee.


\end{abstract}



\section{Introduction}

The purposes of this paper are giving bijective  proofs of certain identities of Andrews-Yee in  \cite{ay}
and extending their results on the  generalized bivariate mock theta functions $\omega(z;q)$ and $\nu(z;q)$ to
the trivariate generalizations $\omega(y,z;q)$ and $\nu(y,z;q)$.
To this end, we first introduce some  definitions and notation.

A \emph{partition} \cite{gea1998} of $n$ is a finite nonincreasing sequence of positive integers
$(\lambda_1, \lambda_2, \ldots, \lambda_\ell)$ such that $n=\lambda_1+\lambda_2+\cdots+\lambda_\ell$.
We write
$\lambda=(\lambda_1,\lambda_2,\ldots,\lambda_\ell)$
and call $\lambda_i$'s the \emph{parts} of $\lambda$.
The \emph{size} of $\lambda$ is the sum of all parts, which is denoted by $|\lambda|$,
and the \emph{length} of $\lambda$ is the number of parts,
which is denoted by $\ell(\lambda)$.
The \emph{conjugate} of $\lambda$ is the partition
$\lambda'=(\lambda'_1,\lambda'_2,\ldots,\lambda'_{\lambda_{1}})$,
where
$\lambda'_{i}=|\{\lambda_{j}\colon\,\lambda_{j}\geq i,1\leq j \leq \ell \}|$ for $1\leq i \leq \lambda_{1}$.
We say $\lambda$  is a \emph{distinct} partition if $\lambda_1>\lambda_2>\cdots>\lambda_\ell$.
The \emph{Ferrers graph} of $\lambda$ is defined to be a left-justified arrangement of $n$ dots in $\ell$ rows consisting of $\lambda_1,\lambda_2,\ldots,\lambda_\ell$ dots,
thus we graphically say that the number of rows of $\lambda$ is $\ell$ and the number of columns of $\lambda$ is $\lambda_1$.

The \emph{$q$-series} \cite{Gasper} notation is defined as :
\begin{align*}
(a;q)_0&:=1,\\
(a;q)_n&:=(1-a)(1-aq)\cdots (1-aq^{n-1}) \quad \text{for $n\ge 1$},\\
(a;q)_{\infty}&:=\lim_{n\to \infty} (a;q)_n, \quad |q|<1.
\end{align*}

Watson \cite{Wa0} defined two third order mock theta functions $\omega(q)$ and $\nu(q)$ as:
\begin{align*}
\omega(q)&:=\sum_{n=0}^{\infty} \frac{q^{2n^2+2n}}{(q;q^2)_{n+1}^2},
\quad \nu(q):=\sum_{n=0}^{\infty} \frac{q^{n^2+n}}{(-q;q^2)_{n+1}}.
\end{align*}
These two functions were  also discovered in Ramanujan's lost notebook \cite{ABLN5,Rlnb}.
In \cite{gea5}, Andrews introduced two generalized bivariate mock theta functions $\omega(z;q)$ and $\nu(z;q)$ of the form:
\begin{equation}
\omega(z;q):=\sum_{n=0}^{\infty} \frac{z^n q^{2n^2+2n}}{(q;q^2)_{n+1}(zq;q^2)_{n+1}}, \quad
\nu(z; q):=\sum_{n=0}^{\infty} \frac{ q^{n^2+n}}{(-zq;q^2)_{n+1}}, \label{andrews1}
\end{equation}
which give $\omega(q)$ and $\nu(q)$ by setting $z=1$ in \eqref{andrews1}.
Later on, Andrews \cite{gea6} proved that
\begin{align}
  \omega(z;q)&=\sum_{n=0}^{\infty} \frac{z^n q^{n}}{(q;q^2)_{n+1}},\label{andrews2.1} \\[6pt]
  \nu(z;q)&=\sum_{n=0}^{\infty} (q/z;q^2)_n (-zq)^n,\label{andrews2.2}
\end{align}
and the above results are  extensively studied  by Choi \cite{choi}.

In \cite{ad2}, Andrews, Dixit  and Yee introduced  two partition functions $p_{\omega}(n)$ and $p_{\nu}(n)$
where
$p_{\omega}(n)$ counts the number of partitions of $n$
in which  all odd parts are less than twice the smallest part,
and  $p_{\nu}(n)$ counts the number of distinct partitions of $n$
with the same constraint as $p_\omega(n)$.
They obtained
\begin{align}
\sum_{n=1}^{\infty} p_{\omega}(n)q^n &=q\, \omega(q),\quad
\sum_{n=0}^{\infty} p_{\nu}(n)q^n= \nu(-q). \label{omeganu}
\end{align}
By the definitions of $p_{\omega}(n)$ and $p_{\nu}(n)$ and the forms of $\omega(q)$ and $\nu(q)$,
\eqref{omeganu} can be expressed as
\begin{align}
\sum_{n=1}^{\infty} \frac{q^n}{(q^{n};q)_{n+1} (q^{2n+2};q^2)_{\infty}}&=\sum_{n=0}^{\infty} \frac{q^{2n^2+2n+1}}{(q;q^2)^2_{n+1}}, \label{omega}\\[6pt]
\sum_{n=0}^{\infty} q^n (-q^{n+1};q)_n (-q^{2n+2};q^2)_{\infty}&=\sum_{n=0}^{\infty} \frac{q^{n^2+n}}{(q;q^2)_{n+1}}. \label{nu}
\end{align}
Andrews also mentioned some other combinatorial interpretations of $p_\nu(n)$ in \cite{ad07,adw, adw2}.

Recently, by using analytic  method with $q$-series,
Andrews and Yee  \cite{ay} provided two bivariate generalizations of \eqref{omega} and \eqref{nu} given in the following theorem.

\begin{thm}{\cite[Theorem 1.]{ay}}\label{thm1}
We have
\begin{align}
\sum_{n=1}^{\infty} \frac{q^n}{(zq^n;q)_{n+1} (zq^{2n+2};q^2)_{\infty}}
&=\sum_{n=0}^{\infty} \frac{z^nq^{2n^2+2n+1}}{(q;q^2)_{n+1} (zq;q^2)_{n+1}},\label{eqthm1}
\\[6pt]
\sum_{n=0}^{\infty} q^n(-zq^{n+1};q)_{n} (-zq^{2n+2};q^2)_{\infty} &=\sum_{n=0}^{\infty} \frac{z^n q^{n^2+n}}{(q;q^2)_{n+1}}.\label{eqthm2}
\end{align}
\end{thm}

It is evident that \eqref{omega} and \eqref{nu} are the special cases of \eqref{eqthm1} and \eqref{eqthm2} when $z=1$, respectively.
Simultaneously, by Theorem \ref{thm1},
they also found another bivariate generalization of $\nu(q)$ as follows:
\begin{equation*}
  \nu_1(z;q):=\sum_{n=0}^{\infty}\frac{z^nq^{n^2+n}}{(-q;q^2)_{n+1}}.
\end{equation*}
Consequently they found a representation of $\nu_1(z;q)$ similar to that of $\nu(z;q)$ given by \eqref{andrews2.2}.

\begin{thm}{\cite[Theorem 2.]{ay}}
We have
  \begin{equation}\label{eqv1}
  \nu_1(z;q)=\sum_{n=0}^{\infty}(zq;q^2)_n(-q)^n.
  \end{equation}
\end{thm}

In the concluding remarks of \cite{ay},  Andrews and Yee  stated that the two identities in Theorem \ref{thm1} cried out for bijective proofs.
They also asserted that it is not difficult to prove \eqref{andrews2.1}, \eqref{andrews2.2} and \eqref{eqv1} combinatorially
but without giving the detailed proof.
In this paper, we answer all the above questions.

Furthermore, we may define trivariate generalizations of $\nu(q)$ and $\omega(q)$ as
\begin{equation*}
    \nu(y,z;q):=\sum_{n=0}^{\infty} \frac{y^nz^nq^{n^2+n}}{(yq;q^2)_{n+1}}
\end{equation*}
and
\begin{equation*}
    \omega(y,z;q):=\sum_{n=0}^{\infty}\frac{y^nz^nq^{2n^2+2n}}{(yq;q^2)_{n+1}(zq;q^2)_{n+1}}.
\end{equation*}

Note that $\nu(-z,-z^{-1};q)=\nu(z;q)$, $\nu(-1,-z;q)=\nu_1(z;q)$
and $\omega(1,z;q)=\omega(z;q)$,
which imply that $\nu(y,z;q)$ and $\omega(y,z;q)$
generalize $\nu_1(z;q)$, $\nu(z;q)$ and $\omega(z;q)$, respectively.
Additionally,
$\nu(y,z;q)$ and $\omega(y,z;q)$ can be written in terms of Choi's \cite{choi} functions:
\begin{equation*}
  \bar{\nu}(\alpha,z;q):=\sum_{n=0}^{\infty}\frac{q^{n(n-1)}z^{2n}}{(-\alpha^2z^2/q^3;q^2)_{n+1}}
\end{equation*}
and
\begin{equation*}
  \bar{\omega}(\alpha,z;q):=\sum_{n=0}^{\infty}\frac{q^{2(n-1)^2-6}\alpha^{2n}z^{4(n+1)}}{(z^2/q;q^2)_{n+1}(\alpha^2z^2/q^3;q^2)_{n+1}}.
\end{equation*}
Here for  accuracy, we use $\bar{\nu}(\alpha,z;q)$ and $\bar{\omega}(\alpha,z;q)$ to substitute the original notations $\nu(\alpha,z;q)$ and $\omega(\alpha,z;q)$ appearing in \cite{choi}.
In particular, we have
\begin{equation}\label{nutobarnu}
  \nu(y,z;q)=\bar{\nu}(iq/\sqrt{z},\sqrt{yz}q;q)
\end{equation}
and
\begin{equation}\label{omegatobaromega}
  \omega(y,z;q)=z^{-2}\bar{\omega}(\sqrt{y}q/\sqrt{z},\sqrt{z}q;q).
\end{equation}

By a combinatorial approach,
we find the following uniform representation of $\nu(y,z;q)$.

\begin{thm}\label{thmofnewnu}
We have
  \begin{equation}\label{newnu}
    \nu(y,z;q)=\sum_{n=0}^{\infty}(-zq;q^2)_n(yq)^n.
  \end{equation}
\end{thm}

Note that by setting $z\rightarrow -z^{-1}$ and $y\rightarrow -z$ in \eqref{newnu},
we obtain the combinatorial interpretation for \eqref{andrews2.2},
and by setting $z\rightarrow -z$ and $y\rightarrow -1$ in \eqref{newnu},
we can also explain  \eqref{eqv1} combinatorially.
Furthermore, the right side of \eqref{newnu} equals
$1+(1+z^{-1}q^{-1})\nu_3(\sqrt{y}q,1/\sqrt{yz};q)$,
where
\begin{equation*}
  \nu_3(\alpha,z;q):=\frac{1}{1+\alpha^2z^2/q^3}\sum_{n=1}^{\infty}\frac{\alpha^{2n}}{q^n}(-q^3/(\alpha z)^2;q^2)_n
\end{equation*}
is defined by Choi \cite{choi}.
Hence, with \eqref{nutobarnu}, Theorem \ref{thmofnewnu} also establishes the connection between
$\bar{\nu}(\alpha,z;q)$ and $\nu_3(\alpha,z;q)$.

Reminiscent to $\nu(y,z;q)$, $\omega(y,z;q)$ also has the following identity.

\begin{thm}\label{thmofnewomega}
  We have
  \begin{equation}\label{newomega}
    \omega(y,z;q)=\sum_{n=0}^{\infty}\frac{y^nq^n}{(zq;q^2)_{n+1}}=\sum_{n=0}^{\infty}\frac{z^nq^n}{(yq;q^2)_{n+1}}.
  \end{equation}
\end{thm}

Thus letting $y\rightarrow 1$ in the first equation of \eqref{newomega} gives the combinatorial interpretation for \eqref{andrews2.1}.

For the bijective proof of Theorem \ref{thm1},
notice that dividing $q$ from the both sides of \eqref{eqthm1}, we deduce
\begin{equation}\label{thmeq3}
\sum_{n=0}^{\infty} \frac{q^{n}}{(zq^{n+1};q)_{n+2} (zq^{2n+4};q^2)_{\infty}}
=\sum_{n=0}^{\infty}\frac{z^nq^{2n^2+2n}}{(q;q^2)_{n+1}(zq;q^2)_{n+1}}.
\end{equation}
Thus, in order to prove  \eqref{eqthm1} bijectively,
it is equivalent to find bijective proof of \eqref{thmeq3}.
Based on the results of Theorem \ref{thmofnewnu} and Theorem  \ref{thmofnewomega},
we can easily interpret the right sides of \eqref{eqthm2} and \eqref{thmeq3} combinatorially,
which leads us to the bijective proofs.


The rest of this paper is organized as follows.
By utilizing a variation of Ferrers graphs called odd Ferrers graphs,
the generalized trivariate mock theta functions $\nu(y,z;q)$ and $\omega(y,z;q)$ as well as
the combinatorial proofs of Theorem \ref{thmofnewnu} and Theorem \ref{thmofnewomega}
are given in Section  \ref{newgennu},
then the corollaries of Theorem \ref{thmofnewnu} and Theorem  \ref{thmofnewomega}
give combinatorial interpretations of \eqref{andrews2.1}, \eqref{andrews2.2}
and \eqref{eqv1} and some other identities involving bivariate generalizations of $\omega(q)$ and $\nu(q)$.
In Section \ref{bjisect},
we present our bijective proofs of \eqref{thmeq3} and \eqref{eqthm2}
by two algorithms and list some examples.
Finally, we conclude the paper with some further remarks in Section \ref{secofrem}.

\section{Trivariate generalizations of the mock theta functions $\omega(q)$ and $\nu(q)$}\label{newgennu}
\setcounter{table}{0}
\setcounter{figure}{0}

In this section,
we recall a variation of Ferrers graphs of partitions called odd Ferrers graphs,
which was previously utilized by Berndt-Yee in \cite{by}
and by Andrews in \cite{ad07,adw,adw2}.
In the spirit of double counting on the odd Ferrers graphs with distinct partitions,
we first give the combinatorial proof for Theorem \ref{thmofnewnu}.
The  corollary of Theorem \ref{thmofnewnu} gives both
the necessary preparation for the proof of Theorem \ref{thmofnewomega}
and the combinatorial interpretation of the right side of \eqref{eqthm2}.
Then by a similar analysis on the odd Ferrers graphs with ordinary partitions,
we constructively prove Theorem \ref{thmofnewomega}
whose corollary leads us to
the combinatorial interpretation for the
right side of \eqref{thmeq3}.


Given a partition $\lambda=(\lambda_1,\lambda_2,\ldots,\lambda_\ell)$,
we draw the Ferrers graph of $\lambda$ and replace each dot by a  box,
then put $0$ into the upper left corner box,
1's into the rest of the boxes in either the first column or the first row,
and 2's into each box except for the first column and the first row.
Since the sum of numbers in each row except for the first row is odd ,
we call it the \emph{odd Ferrers graph}  of shape $\lambda$ and size $n$,
where $n$ is the sum of all numbers in the boxes.
Note that the shape completely determines the size of odd Ferrers graph,
 we can use  the shape $\lambda$ to denote the corresponding  odd Ferrers graph by $\f_\lambda$.
Thus the graphic parameters of $\f_\lambda$ are exactly the same as $\lambda$, that is,
the number of rows (or length) $\ell(\f_\lambda)$ of $\f_\lambda$ is $\ell(\lambda)$
and the number of columns of $\f_\lambda$ is $\lambda_1$.
An odd Ferrers graph $\f_\lambda$ is \emph{distinct} if $\lambda$ is a distinct partition.
For example, the odd Ferrers graph below is $\f_{(6,6,3,2)}$ of size 24.

\begin{center}
 \begin{ytableau}
     0 & 1 & 1 & 1 & 1 & 1 \\
     1 & 2 & 2 & 2 & 2 & 2\\
     1 & 2 & 2 \\
     1 & 2
 \end{ytableau}
\end{center}

\subsection{Generalized trivariate  mock theta function $\nu(y,z;q)$}
Denote by $\mathcal{B}_\nu$ the sets of all distinct odd Ferrers graphs.
Particularly,
let $\mathcal{B}_{\nu}(\ell,m,n)$ denote the set of distinct odd Ferrers graphs of size $n$
with $\ell+1$ rows and $m+1$ columns.
Let $b_{\nu}(\ell,m,n)=|\mathcal{B}_{\nu}(\ell,m,n)|$.
By decomposing $\f_\lambda\in\mathcal{B}_{\nu}(\ell,m,n)$
in two different ways,
we give the  proof of Theorem \ref{thmofnewnu}.

\noindent{\emph{Proof of Theorem \ref{thmofnewnu}.}}
We proceed our proof by showing that the both sides of \eqref{newnu} are generating functions of $b_\nu(\ell,m,n)$, i.e.,
  \begin{equation}\label{dprow}
    \sum_{n=0}^{\infty}\frac{y^nz^nq^{n^2+n}}{(yq;q^2)_{n+1}}
    =\sum_{n=0}^{\infty}\sum_{m=0}^{\infty}\sum_{\ell=0}^{\infty}b_\nu(\ell,m,n)z^\ell y^m q^n
    =\sum_{n=0}^{\infty}(-zq;q^2)_n(yq)^n.
  \end{equation}

\begin{figure}
  \centering
  \includegraphics[width=0.65\textwidth]{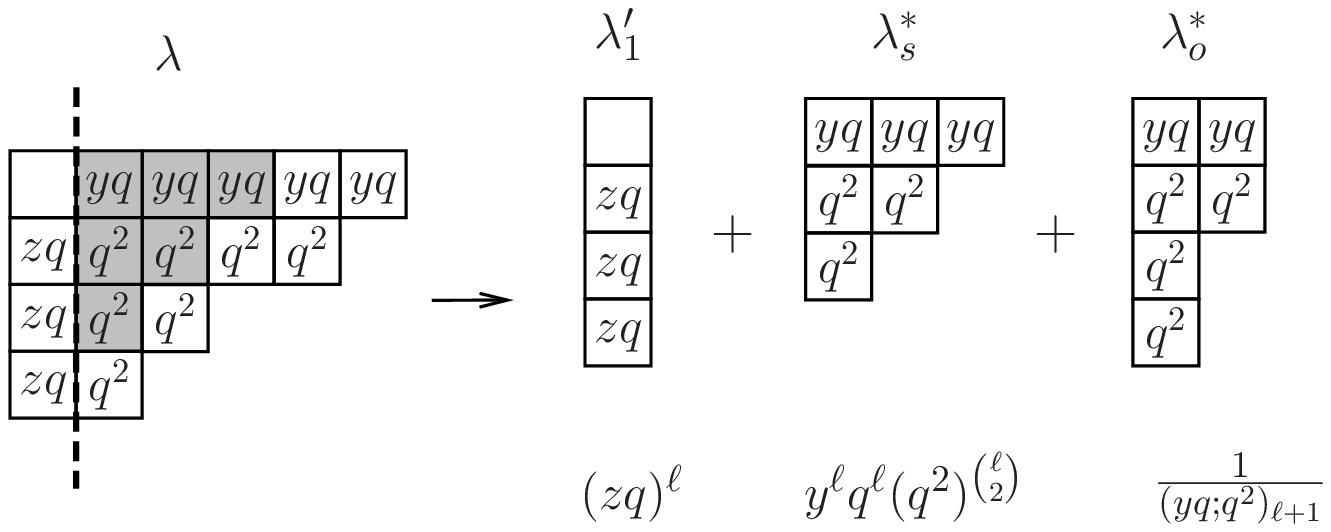}
  \caption{the decomposition of $\f_\lambda$ in the first equation of \eqref{dprow}}
  \label{decompofgfB}
\end{figure}

First, for clarity, we rewrite the left part of \eqref{dprow} as
\begin{equation*}
  \sum_{n=0}^{\infty} \frac{y^nz^n q^{n^2+n}}{(yq;q^2)_{n+1}}=\sum_{\ell=0}^\infty\frac{(zq)^\ell (yq)^\ell(q^2)^{\binom{\ell}{2}}}{(yq;q^2)_{\ell+1}}.
\end{equation*}
Given a distinct odd Ferrers graph $\f_\lambda=\f_{(\lambda_1,\lambda_2,\ldots,\lambda_{\ell+1})}$,
we can decompose $\f_\lambda$ as illustrated in Figure \ref{decompofgfB}.
The first column of $\f_\lambda$ is generated by $(zq)^\ell$,
where the power of $z$ represents the number of rows of $\f_\lambda$ minus 1.
Then by deleting the first column, i.e., the boxes in the $\lambda'_1$,
we obtain a new variational Ferrers graph defined by the distinct partition $\lambda^*=(\lambda_1^*,\lambda_2^*,\ldots,\lambda_p^*)$,
where the boxes in $\lambda_1^*$ are filled by 1's and the boxes in $\lambda_2^*,\ldots,\lambda_p^*$ are filled by 2's.
Furthermore, since $\lambda_1>\lambda_2>\cdots>\lambda_{\ell+1}\geq1$,
it follows that $\ell \leq p \leq \ell+1$ and $\lambda_i^*=\lambda_i-1$ for $1\leq i\leq p$.
Hence by taking out $\ell+1-i$ boxes from $\lambda^*_i$ for $1\leq i\leq \ell$,
we can split $\lambda^*$ into two partitions $\lambda_s^*$ and $\lambda_o^*$,
where $\lambda_s^*=(\ell,\ell-1,\ldots,1)$ is a staircase partition,
and $\lambda_o^*$ is an ordinary partition with no more than $\ell+1$ parts.
Note that each box in the first rows of $\lambda_s^*$ and $\lambda_o^*$ are marked by $yq$ and rest of the boxes of $\lambda_d^*$ and $\lambda_o^*$ are marked by $q^2$'s.
Thus,  the generating function of the variational Ferrers graph of shape $\lambda_s^*$ is $$(yq)^\ell(q^2)^{\sum_{i=1}^{\ell-1}2i}=(yq)^\ell(q^2)^{\binom{\ell}{2}}.$$
Since the number of rows in the variational Ferrers graph of shape $\lambda_o^*$ does not exceed $\ell+1$,
the generating function  is
$$\frac{1}{1-yq}\cdot\frac{1}{(1-yq\cdot q^2)}\cdots\frac{1}{(1-yq\cdot (q^2)^\ell)}=\frac{1}{(yq;q^2)_{\ell+1}}.$$
It is obvious that the power of $y$ indeed equals $\lambda_1-1$, the number of columns of $\f_\lambda$ minus 1.
Hence, combining the three generating functions above, we see that first equation of \eqref{dprow} holds.

\begin{figure}
  \centering
  \includegraphics[width=0.65\textwidth]{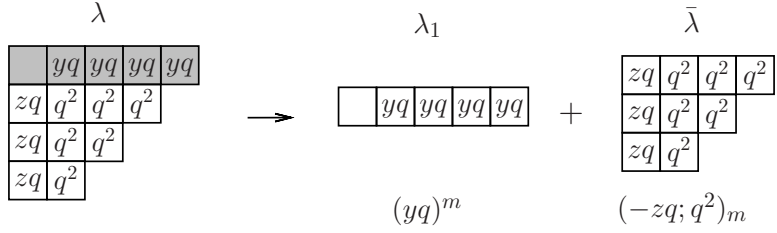}
  \caption{the decomposition of $\f_\lambda$ in the second equation of \eqref{dprow}}
  \label{decompofdprow}
\end{figure}

To complete the proof, it  remains to show that
\begin{equation}\label{gfB2}
  \sum_{m=0}^{\infty}(-zq;q^2)_m (yq)^m=\sum_{n=0}^{\infty}\sum_{m=0}^{\infty}\sum_{\ell=0}^{\infty}b_{\nu}(\ell,m,n)z^\ell y^m q^n.
\end{equation}
Let $\f_\lambda=\f_{(\lambda_1,\lambda_2,\ldots,\lambda_p)}$ be a distinct odd Ferrers graph with $\lambda_1=m+1$.
As shown in Figure \ref{decompofdprow},
all boxes but the first one in $\lambda_1$ are marked by $yq$'s,
implying that the first row of $\f_\lambda$ is  generated by $(yq)^m$.
By deleting the first row, consider the variational Ferrers graph with shape $\bar{\lambda}=(\lambda_2,\lambda_3,\ldots,\lambda_p)$,
where the first box of each row is filled by 1 and the rest of the boxes are filled by 2's.
Thus, we can use $zq$ to mark the first box of each row, and $q^2$'s to mark the other boxes,
which means that  the power of $z$ is the  number of rows of $\f_\lambda$ minus 1.
Since we have $m\geq\lambda_2>\lambda_3>\cdots>\lambda_p\geq1$,
the  variational Ferrers graph of shape $\bar{\lambda}$ is generated by
$$(1+zq)(1+zq\cdot q^2)\cdots(1+zq\cdot (q^2)^{m-1})=(-zq;q^2)_m.$$
Thus, by combining the generating functions above, we complete the proof of \eqref{gfB2} then \eqref{dprow}.

\qed

By specifying the values of $y$ and $z$ in \eqref{newnu},
we can obtain some similar identities related to the
bivariate generalizations of the mock theta function $\nu(q)$.
Additionally, these identities are of their own combinatorial significance.


\begin{cor}\label{coroofdp}
  We have
  \begin{align}
    \sum_{n=0}^{\infty} \frac{z^nq^{n^2+n}}{(q;q^2)_{n+1}}&=\sum_{n=0}^{\infty}(-zq;q^2)_n q^n,\label{rightofnu}\\[6pt]
    \sum_{n=0}^{\infty} \frac{z^nq^{n^2+n}}{(zq;q^2)_{n+1}}&=\sum_{n=0}^{\infty}(-q;q^2)_n(zq)^n,\label{dpcolumn} \\[6pt]
    \sum_{n=0}^{\infty}\frac{q^{n^2+n}}{(zq;q^2)_{n+1}}&=\sum_{n=0}^{\infty}(-q/z;q^2)_n(zq)^n,\label{dpcolumn-row}\\[6pt]
    \sum_{n=0}^{\infty} \frac{z^{2n} q^{n^2+n}}{(zq;q^2)_{n+1}}&=\sum_{n=0}^{\infty}(-zq;q^2)_n(zq)^n. \label{dpcolumn+row}
  \end{align}
\end{cor}
\pf
Setting $y\rightarrow1$ in \eqref{newnu} gives \eqref{rightofnu} and
setting $z\rightarrow 1$, $y\rightarrow z$ in \eqref{newnu} gives \eqref{dpcolumn}.
Let $z\rightarrow z^{-1}$ and $y\rightarrow z$ in \eqref{newnu}, then we see \eqref{dpcolumn-row}.
Finally by letting $y\rightarrow z$ in \eqref{newnu}, we obtain \eqref{dpcolumn+row}.
\qed

For the combinatorial significance,
referring to the proof of Theorem \ref{thmofnewnu},
we know that all four identities in Corollary \ref{coroofdp} are  refinements of the generating function of distinct odd Ferrers graphs $\f_\lambda$,
where the powers of $q$  record the size of $\f_\lambda$,
but the statistics recorded by the powers of $z$ are different.

Denote by $\mathcal{B}_\nu^1(m,n)$ the set of distinct odd Ferrers graphs  of size $n$
with $m+1$ rows
and $b_\nu^1(m,n)$ the cardinality of $\mathcal{B}_\nu^1(m,n)$,
then \eqref{rightofnu} presents the generating function of $b_\nu^1(m,n)$:
\begin{equation}\label{comtorightnu}
  \sum_{n=0}^{\infty} \frac{z^nq^{n^2+n}}{(q;q^2)_{n+1}}=\sum_{n}^{\infty}\sum_{m}^{\infty}b^1_\nu(m,n)z^m q^n=\sum_{n=0}^{\infty}(-zq;q^2)_n q^n,
\end{equation}
which gives the combinatorial interpretation for the right side of \eqref{eqthm2}.
We can also deduce that
\begin{equation}\label{comtonv1}
  \nu_1(z;q)=\sum_{n=0}^{\infty}\sum_{m=0}^{\infty}\sum_{\f_\lambda\in\mathcal{B}_\nu^1(m,n)}(-1)^{\sharp(\f_\lambda)}z^mq^n,
\end{equation}
where $\sharp(\f_\lambda)$ is the number of $1$'s in the odd Ferrers graph $\f_\lambda$.
Hence \eqref{comtonv1} gives the combinatorial interpretation of \eqref{eqv1}.

Let $\mathcal{B}_\nu^2(m,n)$ be the set of distinct odd Ferrers graphs of size $n$
with $m+1$ columns,
and let $b_\nu^2(m,n)=|\mathcal{B}_\nu^2(m,n)|$.
Thus \eqref{dpcolumn} is the generating function of $b_\nu^2(m,n)$:
\begin{equation}\label{comtodpcolumn}
  \sum_{n=0}^{\infty}\frac{z^nq^{n^2+n}}{(zq;q^2)_{n+1}}=\sum_{n=0}^{\infty}\sum_{m=0}^{\infty}b^2_\nu(m,n)z^mq^n=\sum_{n=0}^{\infty}(-q;q^2)_n(zq)^n.
\end{equation}

Let $\mathcal{B}_\nu^3(m,n)$ be the set of distinct odd Ferreres graphs  of size $n$
satisfying that the difference between the number of columns
and the number of rows is $m$.
Denote by $b_\nu^3(m,n)$ the cardinality of $\mathcal{B}_\nu^3(m,n)$,
then \eqref{dpcolumn-row} is the generating function of $b_\nu^3(m,n)$:
\begin{equation}\label{comtodpcolumn-row}
  \sum_{n=0}^{\infty}\frac{q^{n^2+n}}{(zq;q^2)_{n+1}}=\sum_{n=0}^{\infty}\sum_{m=0}^{\infty}b^3_\nu(m,n)z^mq^n=\sum_{n=0}^{\infty}(-q/z;q^2)_n(zq)^n.
\end{equation}
For any $\f_\lambda\in\mathcal{B}_\nu^3(m,n)$, since $\lambda$ is a distinct partition, we have $\lambda_1\geq\ell(\lambda)$ so that the range of $m$ begins with 0.
Similar to \eqref{comtonv1}, in terms of the set $\mathcal{B}_\nu^3(m,n)$, we can explain \eqref{andrews2.2} combinatorially as
\begin{equation}\label{comtonv}
  \nu(z;q)=\sum_{n=0}^{\infty}\sum_{m=0}^{\infty}\sum_{\f_\lambda\in\mathcal{B}_\nu^3(m,n)}(-1)^{\sharp(\f_\lambda)}z^mq^n.
\end{equation}

Let $\mathcal{B}_\nu^4(m,n)$ be the set of distinct odd Ferrers graphs $\f_\lambda$ of size $n$ and $\sharp(\f_\lambda)=m$, and let $b_\nu^4(m,n)=|\mathcal{B}_\nu^4(m,n)|$. By \eqref{dpcolumn+row}, we have
\begin{equation}\label{comtodpcolumn+row}
  \sum_{n=0}^{\infty} \frac{z^{2n} q^{n^2+n}}{(zq;q^2)_{n+1}}=\sum_{n=0}^{\infty}\sum_{m=0}^{\infty}b^4_\nu(m,n)z^mq^n=\sum_{n=0}^{\infty}(-zq;q^2)_n(zq)^n.
\end{equation}

In light of \eqref{comtonv1} and \eqref{comtonv}, it follows from  \eqref{comtodpcolumn} and \eqref{comtodpcolumn+row} that
\begin{equation}\label{comtodpcol-}
  \sum_{n=0}^{\infty}\frac{z^nq^{n^2+n}}{(-zq;q^2)_{n+1}}=\sum_{n=0}^{\infty}\sum_{m=0}^{\infty}\sum_{\f_\lambda\in\mathcal{B}_\nu^2(m,n)}(-1)^{\sharp(\f_\lambda)}z^mq^n=\sum_{n=0}^{\infty}(q;q^2)_n(-zq)^n
\end{equation}
and
\begin{equation}\label{comtodpcol+row-}
  \sum_{n=0}^{\infty} \frac{z^{2n} q^{n^2+n}}{(-zq;q^2)_{n+1}}=\sum_{n=0}^{\infty}\sum_{m=0}^{\infty}(-1)^mb^4_\nu(m,n)z^mq^n=\sum_{n=0}^{\infty}(zq;q^2)_n(-zq)^n,
\end{equation}
where \eqref{comtodpcol-} and \eqref{comtodpcol+row-} give the combinatorial proofs of the identities deduced from \eqref{dpcolumn} and \eqref{dpcolumn+row}
by letting $q\rightarrow -q$, respectively.

%
%

\subsection{Generalized trivariate mock theta function $\omega(y,z;q)$}
In the rest of this section, we study the combinatorial identities related to the trivariate generalization $\omega(y,z;q)$ of $\omega(q)$.
We use $\mathcal{B}_\omega$ to denote
the set of all odd Ferrers graphs.
Let $\mathcal{B}_{\omega}(\ell,m,n)$ denote the set of odd Ferrers graphs of size $n$ consisting of
$\ell+1$ rows and $m+1$ columns,
and let $b_{\omega}(\ell,m,n)=|\mathcal{B}_{\omega}(\ell,m,n)|$.
Now we give the proof of Theorem \ref{thmofnewomega}.

\noindent{\it{Proof of Theorem \ref{thmofnewomega}.}}
By noticing that
\begin{equation}\label{newomega2}
  \sum_{n=0}^{\infty}\frac{y^nz^nq^{2n^2+2n}}{(yq;q^2)_{n+1}(zq;q^2)_{n+1}}=\sum_{n=0}^{\infty}\frac{z^nq^n}{(yq;q^2)_{n+1}}
\end{equation}
is immediately from
\begin{equation}\label{newomega1}
  \sum_{n=0}^{\infty}\frac{y^nz^nq^{2n^2+2n}}{(yq;q^2)_{n+1}(zq;q^2)_{n+1}}=\sum_{n=0}^{\infty}\frac{y^nq^n}{(zq;q^2)_{n+1}}
\end{equation}
since the variables $y$ and $z$ in the
term $\frac{y^nz^nq^{2n^2+2n}}{(yq;q^2)_{n+1}(zq;q^2)_{n+1}}$ are symmetric,
we only prove \eqref{newomega1} by showing that
\begin{equation}\label{oprow}
  \sum_{m=0}^{\infty}\frac{y^mq^{m^2+m}}{(yq;q^2)_{m+1}}\cdot\frac{z^mq^{m^2+m}}{(zq;q^2)_{m+1}}=\sum_{n=0}^{\infty}\sum_{m=0}^{\infty}\sum_{\ell=0}^{\infty}b_\omega(\ell,m,n)z^\ell y^m q^n=\sum_{m=0}^{\infty}\frac{y^mq^m}{(zq;q^2)_{m+1}}.
\end{equation}

\begin{figure}
  \centering
  \includegraphics[width=0.65\textwidth]{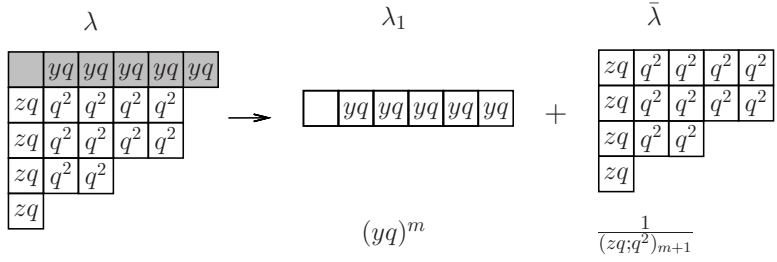}
  \caption{The decomposition of $\f_\lambda$ in the second equation  of \eqref{oprow}}
  \label{decompofoprowright}
\end{figure}

For the second equation of \eqref{oprow},
given an odd Ferrers graph $\f_\lambda=\f_{(\lambda_1,\lambda_2,\ldots,\lambda_p)}$ with $\lambda_1=m+1$,
we can split  $\f_\lambda$ as Figure \ref{decompofoprowright}.
Note that all boxes expect the first one of $\lambda_1$ are filled by 1's
so that we mark each of these box by $yq$,
which implies that $(yq)^m$ generates the nonempty boxes in $\lambda_1$ and
the power of $y$ is the number of columns of $\f_\lambda$ minus 1.
Then let $\bar{\lambda}=(\lambda_2,\ldots,\lambda_p)$ be the shape of the variational Ferrers graph
whose first boxes of all rows are filled by 1's and other boxes are filled by 2's.
Thus, we can use $zq$ to mark the first box in each row and $q^2$'s to mark the other boxes.
It obvious that the power of $q$ is the size of $\f_\lambda$
and the power of $z$ is $p-1$,
which is the number of rows of $\f_\lambda$ minus 1.
Since $m+1\geq\lambda_2\cdots\geq\lambda_p$, the variational Ferrers graph of shape $\bar{\lambda}=(\lambda_2,\ldots,\lambda_p)$
is generated by
$$\frac{1}{(1-zq)}\cdot\frac{1}{(1-zq\cdot q^2)}\cdots\frac{1}{(1-zq\cdot (q^2)^m)}=\frac{1}{(zq;q^2)_{m+1}},$$
implying that the right side of \eqref{newomega1} generates $b_\omega(\ell,m,n)$.

Before proving the first equation of \eqref{oprow},
recall the Frobenius symbol of partitions.
Let $\lambda=(\lambda_1,\lambda_2,\ldots,\lambda_p)$ be a partition
whose Durfee square has side length $m$,
where the \emph{Durfee square} is the largest square contained in the Ferrers graph of $\lambda$.
The \emph{Frobenius symbol} of $\lambda$ is a two-row array of distinct partitions as the following form
\begin{equation*}
  \left(\mu_1,\mu_2,\ldots,\mu_m \atop \nu_1,\nu_2,\ldots,\nu_m\right),
\end{equation*}
where $\mu_i=\lambda_i-i$ and $\nu_i=\lambda'_i-i$ for $1\leq i\leq m$.
Thus we have $\mu_1>\mu_2>\cdots>\mu_m\geq0$ and $\nu_1>\nu_2>\cdots>\nu_m\geq0$.
Notice that the partition is uniquely determined by its Frobenius symbol.

\begin{figure}
  \centering
  \includegraphics[width=0.65\textwidth]{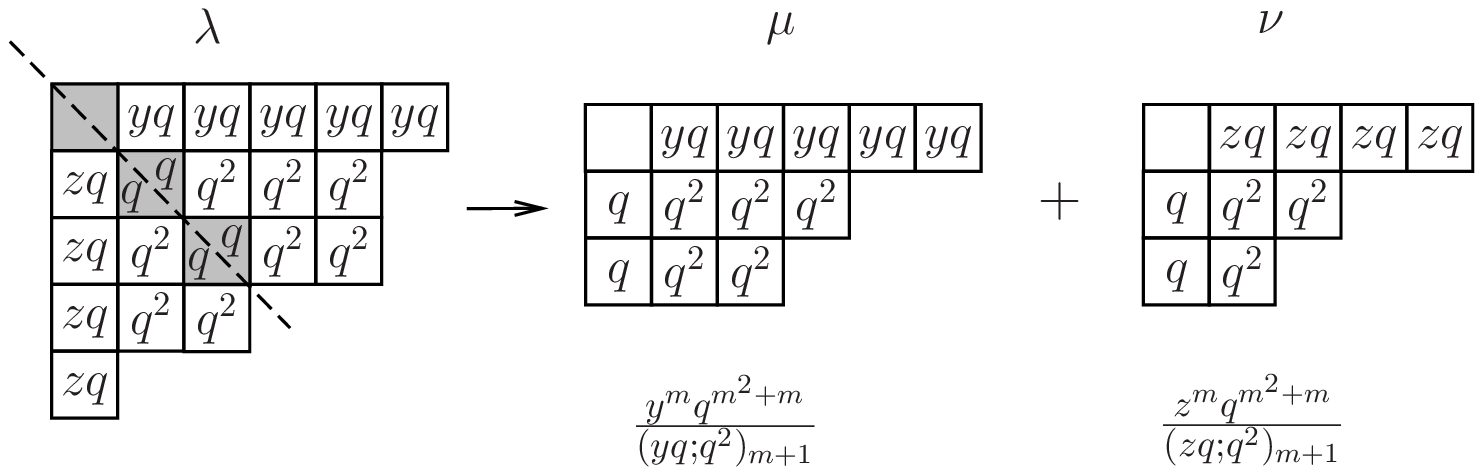}
  \caption{the decomposition of $\f_\lambda$ in the first equation of \eqref{oprow}}
  \label{decompofoprowleft}
\end{figure}

In terms of the Frobenius symbol,
given an odd Ferrers graph
$\f_\lambda=\f_{(\lambda_1,\lambda_2,\ldots,\lambda_p)}$ of shape $\lambda$
whose  Durfee square has side length  $m+1$,
we decompose $\f_\lambda$ as Figure \ref{decompofoprowleft}.
To be more specific, we split $\f_\lambda$
by the diagonal of the Durfee square into
two odd Ferrers graphs $\f_\mu$ and $\f_\nu$.
Note that $\mu=(\mu_1,\mu_2,\ldots,\mu_{m+1})$ and $\nu=(\nu_1,\nu_2,\ldots,\nu_{m+1})$ are both distinct partitions with the same length $m+1$,
where $\mu_i=\lambda_i-i+1$ and $\nu_i=\lambda'_i-i+1$ for $1\leq i\leq m+1$.
Thus $\f_\mu$ and $\f_\nu$ are both distinct.
Note that except for the first box, all boxes in the first row of $\f_\lambda$ are marked by $yq$'s,
then it follows that the power of $y$ is $\lambda_1-1$, which is the number of columns of $\f_\mu$ minus 1.
Hence by \eqref{comtodpcolumn}, $\f_\mu$ is generated by $\frac{y^mq^{m^2+m}}{(yq;q^2)_{m+1}}$.
Furthermore, since $zq$ is used to mark  all the boxes in the first column of $\f_\lambda$
except for the first one,
the power of $z$ now is records the number of the columns of $\f_\nu$ minus 1,
which implies that  $\f_\nu$ is generated by
$\frac{z^mq^{m^2+m}}{(zq;q^2)_{m+1}}$ also by \eqref{comtodpcolumn}.
Therefore, by combining these two above generating functions,
we have shown that the left side of \eqref{newomega1} also generates $b_{\omega}(\ell,m,n)$,
which implies that \eqref{oprow} holds.
\qed

Similar to Corollary \ref{coroofdp} of Theorem \ref{thmofnewnu},
for the ordinary partition case in Theorem \ref{thmofnewomega},
we also derive the following identities concerned with the bivariate generalizations of the mock theta function $\omega(q)$.

\begin{cor}\label{coroofop}
  We have
  \begin{align}
    \sum_{n=0}^{\infty}\frac{z^nq^{2n^2+2n}}{(q;q^2)_{n+1}(zq;q^2)_{n+1}}&=\sum_{n=0}^{\infty}\frac{q^n}{(zq;q^2)_{n+1}}=\sum_{n=0}^{\infty}\frac{z^nq^n}{(q;q^2)_{n+1}}, \label{opcolumn} \\[6pt]
    \sum_{n=0}^{\infty}\frac{q^{2n^2+2n}}{(q/z;q^2)_{n+1}(zq;q^2)_{n+1}}&=\sum_{n=0}^{\infty}\frac{z^{-n}q^n}{(zq;q^2)_{n+1}}=\sum_{n=0}^{\infty}\frac{z^nq^n}{(q/z;q^2)_{n+1}}, \label{opcolumn-row} \\[6pt]
    \sum_{n=0}^{\infty}\left(\frac{z^nq^{n^2+n}}{(zq;q^2)_{n+1}}\right)^2&=\sum_{n=0}^{\infty}\frac{z^nq^n}{(zq;q^2)_{n+1}}. \label{opcolumn+row}
  \end{align}
\end{cor}

\pf Let $y\rightarrow 1$, $y\rightarrow z^{-1}$ and $y\rightarrow z$ in \eqref{newomega}, respectively,
then we obtain \eqref{opcolumn}, \eqref{opcolumn-row} and \eqref{opcolumn+row}.
\qed

Analogous to the proof of Theorem \ref{thmofnewomega},
the above identities are refinements of the generating functions of odd Ferrers graphs $\f_\lambda$.
The power of $q$ in each equation represents the size of $\f_\lambda$
but the power of $z$ records different statistics of $\f_\lambda$.

Let $\mathcal{B}_\omega^1(m,n)$ be the set of
odd Ferrers graphs of size $n$ with $m+1$ rows
and $\mathcal{B}_\omega^{1'}(m,n)$ be the set of
odd Ferrers graphs of size $n$ with $m+1$ columns.
Denote by $b_\omega^1(m,n)$ and $b_\omega^{1'}(m,n)$
the cardinalities of $\mathcal{B}_\omega^1(m,n)$ and $\mathcal{B}_\omega^{1'}(m,n)$, respectively.
Then the identities in \eqref{opcolumn} can be
 written as
\begin{equation}\label{comtooprow}
  \sum_{n=0}^{\infty}\frac{z^nq^{2n^2+2n}}{(q;q^2)_{n+1}(zq;q^2)_{n+1}}
   =\sum_{n=0}^{\infty}\sum_{m=0}^{\infty}b_\omega^1(m,n)z^m q^n
   =\sum_{n=0}^{\infty}\frac{q^n}{(zq;q^2)_{n+1}}
\end{equation}
and
\begin{equation}\label{comtoopcol}
  \sum_{n=0}^{\infty}\frac{z^nq^{2n^2+2n}}{(q;q^2)_{n+1}(zq;q^2)_{n+1}}
   =\sum_{n=0}^{\infty}\sum_{m=0}^{\infty}b_\omega^{1'}(m,n)z^m q^n
   =\sum_{n=0}^{\infty}\frac{z^nq^n}{(q;q^2)_{n+1}},
\end{equation}
where \eqref{comtooprow} can be  utilized  to explain the right side of \eqref{thmeq3}
and \eqref{comtoopcol} gives the combinatorial interpretation for \eqref{andrews2.1}.
Note here we can directly deduce
\begin{equation*}
  \sum_{n=0}^{\infty}\frac{q^n}{(zq;q^2)_{n+1}}=\sum_{n=0}^{\infty}\frac{z^nq^n}{(q;q^2)_{n+1}}
\end{equation*}
without the help of Theorem \ref{thmofnewomega}
since conjugating partitions gives
a bijection between $\mathcal{B}_\omega^1(m,n)$ and $\mathcal{B}_\omega^{1'}(m,n)$.

Let $\mathcal{B}_\omega^2(m,n)$ be the set of odd Ferrers graphs of size $n$ whose
difference between the number of rows
and the number of columns is $m$,
and let $\mathcal{B}_\omega^{2'}(m,n)$ be the set of odd Ferrers graphs of size $n$ whose
difference between the number of columns
and the number of rows is $m$.
Setting $b_\omega^2(m,n)=|\mathcal{B}_\omega^2(m,n)|$ and $b_\omega^{2'}(m,n)=|\mathcal{B}_\omega^{2'}(m,n)|$,
gives that \eqref{opcolumn-row} are the generating functions of $b_\omega^2(m,n)$ and $b_\omega^{2'}(m,n)$:
\begin{equation}\label{comtooprow-col}
  \sum_{n=0}^{\infty}\frac{q^{2n^2+2n}}{(q/z;q^2)_{n+1}(zq;q^2)_{n+1}}
  =\sum_{n=0}^{\infty}\sum_{m=-\infty}^{\infty}b^2_\omega(m,n)z^mq^n
  =\sum_{n=0}^{\infty}\frac{z^{-n}q^n}{(zq;q^2)_{n+1}}
\end{equation}
and
\begin{equation}\label{comtoopcol-row}
  \sum_{n=0}^{\infty}\frac{q^{2n^2+2n}}{(q/z;q^2)_{n+1}(zq;q^2)_{n+1}}
  =\sum_{n=0}^{\infty}\sum_{m=-\infty}^{\infty}b^{2'}_\omega(m,n)z^mq^n
  =\sum_{n=0}^{\infty}\frac{z^nq^n}{(q/z;q^2)_{n+1}}.
\end{equation}
Taking the conjugate of a partition also gives a bijection between $\mathcal{B}_\omega^2(m,n)$
and $\mathcal{B}_\omega^{2'}(m,n)$,
implying that we can obtain the following identity in a purely combinatorial way:
\begin{equation*}
  \sum_{n=0}^{\infty}\frac{z^{-n}q^n}{(zq;q^2)_{n+1}}=\sum_{n=0}^{\infty}\frac{z^nq^n}{(q/z;q^2)_{n+1}}.
\end{equation*}

Let $\mathcal{B}_\omega^3(m,n)$ be the set of odd Ferrers graphs $\f_\lambda$ of size $n$ and $\sharp(\f_\lambda)=m$,
then the combinatorial proof of \eqref{opcolumn+row} is presented by
\begin{equation}\label{comtooprow+col}
  \sum_{n=0}^{\infty}\left(\frac{z^nq^{n^2+n}}{(zq;q^2)_{n+1}}\right)^2=\sum_{n=0}^{\infty}\sum_{m=0}^{\infty}b_\omega^3(m,n)z^mq^n=\sum_{n=0}^{\infty}\frac{z^nq^n}{(zq;q^2)_{n+1}},
\end{equation}
where $b_\omega^3(m,n)=|\mathcal{B}_\omega^3(m,n)|$.

By either letting $q\rightarrow-q$ in Corollary \ref{coroofop}
or analyzing \eqref{comtooprow}--\eqref{comtooprow+col} combinatorially,
we can easily arrive at the following corollary.
\begin{cor}
   We have
  \begin{align*}
    \sum_{n=0}^{\infty}\frac{z^nq^{2n^2+2n}}{(-q;q^2)_{n+1}(-zq;q^2)_{n+1}}&=\sum_{n=0}^{\infty}\sum_{m=0}^{\infty}\sum_{\f_\lambda\in\mathcal{B}_\omega^{1}(m,n)}(-1)^{\sharp(\f_\lambda)}z^m q^n=\sum_{n=0}^{\infty}\frac{(-q)^n}{(-zq;q^2)_{n+1}}, \\[6pt]
    \sum_{n=0}^{\infty}\frac{z^nq^{2n^2+2n}}{(-q;q^2)_{n+1}(-zq;q^2)_{n+1}}&=\sum_{n=0}^{\infty}\sum_{m=0}^{\infty}\sum_{\f_\lambda\in\mathcal{B}_\omega^{1'}(m,n)}(-1)^{\sharp(\f_\lambda)}z^m q^n=\sum_{n=0}^{\infty}\frac{z^n(-q)^n}{(-q;q^2)_{n+1}}, \\[6pt]
    \sum_{n=0}^{\infty}\frac{q^{2n^2+2n}}{(-q/z;q^2)_{n+1}(-zq;q^2)_{n+1}}&=\sum_{n=0}^{\infty}\sum_{m=-\infty}^{\infty}\sum_{\f_\lambda\in\mathcal{B}_\omega^{2}(m,n)}(-1)^{\sharp(\f_\lambda)}z^m q^n=\sum_{n=0}^{\infty}\frac{z^{-n}(-q)^n}{(-zq;q^2)_{n+1}}, \\[6pt]
    \sum_{n=0}^{\infty}\frac{q^{2n^2+2n}}{(-q/z;q^2)_{n+1}(-zq;q^2)_{n+1}}&=\sum_{n=0}^{\infty}\sum_{m=-\infty}^{\infty}\sum_{\f_\lambda\in\mathcal{B}_\omega^{2'}(m,n)}(-1)^{\sharp(\f_\lambda)}z^m q^n=\sum_{n=0}^{\infty}\frac{z^n(-q)^n}{(-q/z;q^2)_{n+1}}, \\[6pt]
    \sum_{n=0}^{\infty}\left(\frac{z^nq^{n^2+n}}{(-zq;q^2)_{n+1}}\right)^2&=\sum_{n=0}^{\infty}\sum_{m=0}^{\infty}(-1)^mb_\omega^3(m,n)z^mq^n=\sum_{n=0}^{\infty}\frac{z^n(-q)^n}{(-zq;q^2)_{n+1}}.
  \end{align*}
\end{cor}

\begin{remark}
The three identities \eqref{andrews2.1}, \eqref{andrews2.2} and \eqref{eqv1} of
$\omega(z;q)$, $\nu(z;q)$ and $\nu_1(z;q)$ are mentioned in \cite[Eq. (2) and (9)]{ay}
but with no combinatorial proof.
Later Chern \cite{chern} proved \eqref{andrews2.1}, \eqref{andrews2.2} and \eqref{eqv1}
by using two completely different bijections.
Here by utilizing the method of double counting on only one combinatorial structure,
the odd Ferrers graph,
we not only unify the proofs of  these three identities given by \eqref{comtoopcol}, \eqref{comtonv} and \eqref{comtonv1}, respectively,
but also obtain several other identities
involving  the  bivariate
generalizations of the mock theta functions $\omega(q)$ and $\nu(q)$
with their own combinatorial significance.
\end{remark}

\section{Bijective proofs of Theorem \ref{thm1}}\label{bjisect}
\setcounter{table}{0}
\setcounter{figure}{0}

Built on the combinatorial analysis of $\omega(y,z;q)$ and $\nu(y,z;q)$,
 in this section, we  present the bijective proofs for \eqref{thmeq3} and \eqref{eqthm2}.

Let $\mathcal{P}_{\omega}$ be the set of partitions
with unique smallest part which can be 0 and satisfying that
all the odd parts do not exceed twice  the smallest part plus 1.
Let $\mathcal{P}_{\omega }(m,n)$ denote the  set of partitions in $\mathcal{P}_{\omega}$ with length $m+1$ and size $n$,
and $p_{\omega}(m,n)$ be the cardinality of $\mathcal{P}_{\omega}(m,n)$.
Similarly, let $\mathcal{P}_{\nu}$ denote the set of distinct partitions whose smallest part can be 0
satisfying that  all odd parts are less than twice  the smallest part,
and let $\mathcal{P}_{\nu}(m,n)$ be the set of  such partitions with length $m+1$ and size $n$.
Denote by $p_{\nu}(m,n)$ the cardinality of $\mathcal{P}_{\nu}(m,n)$.
By these definitions, it directly follows that
the left sides of \eqref{thmeq3} and \eqref{eqthm2}
are the generating functions of $p_\omega(m,n)$ and $p_\nu(m,n)$, respectively, that is,
\begin{align}
\sum_{n=0}^{\infty} \frac{q^{n}}{(zq^{n+1};q)_{n+2} (zq^{2n+4};q^2)_{\infty}}
&=\sum_{n=0}^{\infty}\sum_{m=0}^{\infty}p_{\omega}(m,n)z^mq^n, \label{pomega}\\[6pt]
\sum_{n=0}^{\infty} q^{n}(-zq^{n+1};q)_{n} (-zq^{2n+2};q^2)_{\infty}
&=\sum_{n=0}^{\infty}\sum_{m=0}^{\infty}p_{\nu}(m,n)z^mq^n.\label{pnu}
\end{align}
Therefore, by \eqref{comtooprow} and \eqref{comtorightnu},
finding a bijective proof of \eqref{thmeq3} is to
establish a one-to-one correspondence between the set $\mathcal{P}_\omega(m,n)$
and the set $\mathcal{B}_\omega^1(m,n)$,
similarly finding a bijective proof of \eqref{eqthm2} is
to establish a one-to-one correspondence between $\mathcal{P}_\nu(m,n)$
and $\mathcal{B}_\nu^1(m,n)$ for all $m,n\geq0$.

Given a partition $\lambda=(\lambda_1,\lambda_2,\ldots,\lambda_\ell)$,
define $\lambda^+=(\lambda_1,\ldots,\lambda_\ell,0)$
and $\lambda^-=(\lambda_1,\ldots,\lambda_{\ell-1})$.
Additionally, we  need to introduce the following operators acting on partition $\lambda$:
\begin{description}
  \item[\bm{$\phi^+$}:] define $\phi^{+}(\lambda)=(\lambda_1,\ldots,\lambda_{\ell-1},\lambda_\ell+1)_{\geq}$,

  \item[\bm{$\phi^-$}:] define $\phi^{-}$($\lambda$)=$(\lambda_1,\ldots,\lambda_{\ell-1},\lambda_{\ell}-1)$,

  \item[\bm{$\phi^{+}_c$}:] suppose that $\lambda_i$ is one of the largest odd parts of $\lambda$,
   define $\phi^{+}_c$($\lambda$)=$(\lambda_1,\ldots,\lambda_{i-1},\lambda_{i+1},\ldots,\\\lambda_{\ell},(\lambda_i+1)/2,(\lambda_i-1)/2)_{\geq}$,

  \item[\bm{$\phi^-_c$}:] define $\phi^-_c$($\lambda$)=$(\lambda_1,\ldots,\lambda_{\ell-2},\lambda_{\ell-1}+\lambda_{\ell})_{\geq}$,

  \item[\bm{$\phi^{+}_e$}:] define $\phi^{+}_e$($\lambda$)=$(\lambda_1+2,\ldots,\lambda_{\ell}+2)$,

  \item[\bm{$\phi^-_e$}:] define $\phi^-_e$($\lambda$)=$(\lambda_1-2,\ldots,\lambda_{\ell}-2)$,

  \item[\bm{$\phi^{+}_o$}:] define $\phi^{+}_o(\lambda)=(\lambda_1+2,\ldots,\lambda_{\ell-1}+2,\lambda_{\ell}+1)$,

  \item[\bm{$\phi^-_o$}:] define $\phi^-_o(\lambda)=(\lambda_1-2,\ldots,\lambda_{\ell-1}-2,\lambda_{\ell}-1)_{\geq}$,

  \item[\bm{$\phi^{*}$}:] define $\phi^{*}(\lambda)=(\lambda_1-1,\ldots,\lambda_{\ell}-1)$,
\end{description}
where the subscript $\geq$ forces the numbers in parentheses rearranged in nonincreasing order.

For any operator $\bullet$ acting on $\lambda$,
let $d_\bullet(\lambda)$ be the difference bewteen the size of the original partition $\lambda$ and the resulting partition $\bullet(\lambda)$,
i.e., $d_\bullet(\lambda)=|\lambda|-|\bullet(\lambda)|$.

\subsection{Bijection between $\mathcal{P}_{\omega }(m,n)$ and $\mathcal{B}_\omega^1(m,n)$}

We begin with providing  some necessary lemmas for the bijective proof of \eqref{thmeq3}.

\begin{lem}\label{lem1}
Given any $\lambda=(\lambda_1,\ldots,\lambda_\ell)\in\mathcal{P}_{\omega}$,
define the \emph{ destructive} operator $\psi^-$ as follows.
\begin{enumerate}[i.]
  \item $\psi^-(\lambda)=\phi^-(\lambda^-)$ if $\lambda_{\ell}=0$ and
      $\lambda_{\ell-1}\geq2$;
  \item $\psi^-(\lambda)=\phi^-(\phi^-_c(\lambda))$ if $\lambda_{\ell}\geq 0$ and $\lambda_{\ell-1}=\lambda_{\ell}+1$;
  \item $\psi^-(\lambda)=\phi^-_o(\lambda)$ if $\lambda_{\ell}\geq1$ and $\lambda_{\ell-1}\geq \lambda_{\ell}+2$.
\end{enumerate}
Then we have that $d_{\psi^-}(\lambda)$ is odd
and $\psi^-(\lambda)\in\mathcal{P}_\omega$.
\end{lem}

\pf
Recall that $\mathcal{P}_{\omega}$ is the  set of  partitions
with only one smallest part,
which can be 0,
satisfying that odd parts are no more than twice the smallest part plus 1, that is,
$\lambda_{\ell-1}>\lambda_\ell\geq 0$
and $\lambda_i\leq 2\lambda_\ell+1$
if $\lambda_i$ is odd for $1\leq i\leq \ell-1$.
Note that the actions of $\phi^-$ and $\phi^-_o$ already guarantee that
$d_{\psi^-}(\lambda)$ is odd and
the property of  unique smallest part
in the resulting partitions.
Thus, depending on the numerical relationship between $\lambda_{\ell-1}$ and $\lambda_\ell$,
we only need to verify the constraint of odd parts by the following three cases.

Assuming that $\lambda_\ell=0$
and $\lambda_{\ell-1}\geq2$, by the definition of $\mathcal{P}_{\omega}$,
all parts in $\lambda$ are even,
since  $\phi^-(\lambda^-)=(\lambda_1,\lambda_2,\ldots,\lambda_{\ell-1}-1)$,
which implies that all parts except for $\lambda_{\ell-1}-1$ in $\phi^-(\lambda^-)$ are even.
Thus, we have $\phi^-(\lambda^-)\in\mathcal{P}_{\omega}$.
For example, let $\lambda=(6,4,2,2,0)\in\mathcal{P}_\omega$,
then $\psi^-(\lambda)=(6,4,2,1)\in\mathcal{P}_\omega$.

If $\lambda_{\ell}\geq 0$ and $\lambda_{\ell-1}=\lambda_{\ell}+1$,
we may assume that $\lambda_{\ell}=a$ and $\lambda_{\ell-1}=a+1$,
implying that $\lambda_{\ell-1}+\lambda_{\ell}=2a+1$ is odd.
For convenience, set $\lambda_0=\infty$,
then there  exists a unique $i_0$, with $1\leq i_0 \leq \ell-1$,
such that
$\lambda_{i_0}\leq 2a+1< \lambda_{i_0-1}$.
Therefore, we see that
$\phi^-(\phi^-_c(\lambda))=(\lambda_1,\ldots,\lambda_{\ell-2},2a)$
if $i_0=\ell-1$
and
$\phi^-(\phi^-_c(\lambda))=(\lambda_1,\ldots,\lambda_{i_0-1},2a+1,\lambda_{i_0},\ldots,\lambda_{\ell-2}-1)$
if $1\leq i_0\leq \ell-2$.
For the first case,
since $\lambda_{\ell-2}>2a+1$ and
no odd parts in $\lambda$ exceed $2a+1$,
it follows that $\lambda_i$ is even for $1\leq i\leq \ell-2$. Thus $\phi^-(\phi^-_c(\lambda))\in\mathcal{P}_\omega$.
For the second case,
since $\lambda_\ell=a$,
$2a+1$ is also the largest odd part in $\phi^-(\phi^-_c(\lambda))$;
on the other hand,
the smallest part in $\phi^-(\phi^-_c(\lambda))$ is $\lambda_{\ell-2}-1$, which satisfies
$\lambda_{\ell-2}-1\geq\lambda_{\ell-1}-1=a$.
Thus, all odd parts in $\phi^-(\phi^-_c(\lambda))$ do not exceed $2(\lambda_{\ell-2}-1)+1$,
implying that $\phi^-(\phi^-_c(\lambda))\in\mathcal{P}_\omega$.
For example, let $\lambda=(10,8,7,7,5,4)\in\mathcal{P}_\omega$,
then $\psi^-(\lambda)=(10,9,8,7,6)\in\mathcal{P}_\omega$.

If $\lambda_{\ell}\geq 1$ and $\lambda_{\ell-1}\geq \lambda_{\ell}+2$,
then we have $\phi^-_o(\lambda)=(\lambda_1-2,\ldots,\lambda_{\ell-1}-2,\lambda_{\ell}-1)$.
Since for all $1\leq i\leq \ell-1$, $\lambda_i\leq 2\lambda_{\ell}+1$ if $\lambda_i$ is odd, it is clear that $\lambda_i-2\leq 2\lambda_{\ell}-1=2(\lambda_{\ell}-1)+1$,
which means that $\phi^-_o(\lambda)\in\mathcal{P}_{\omega }$.
For example, let $\lambda=(8,7,5,5,3)\in\mathcal{P}_\omega$,
then $\psi^-(\lambda)=(6,5,3,3,2)\in\mathcal{P}_\omega$.
\qed

Here we may explain the motivation of naming $\psi^-$ by \emph{destructive}.
Note that when any $\lambda\in\mathcal{P}_\omega$ is acted on by $\psi^-$, the size of the resulting partition $\psi^-(\lambda)$ is decreased
and in two of the three cases,
the length of the resulting partition $\psi^-(\lambda)$ is decreased by 1,
which implies that by iteratively applying the operator $\psi^-$, the partition $\lambda$
can be converted to
a simpler partition $\mu\in\mathcal{P}_\omega$.
Corresponding to the \emph{destructive} operator $\psi^-$,
we may define two \emph{constructive} operators $\psi_1^+$ and $\psi_2^+$
which have similar properties  to $\psi^-$

\begin{lem}\label{lem2}
Given any  $\lambda=(\lambda_1,\ldots,\lambda_\ell)\in\mathcal{P}_{\omega}$,
define the \emph{constructive} operators $\psi_2^+$ by $\psi_2^+(\lambda)=\phi^+_o(\lambda)$ and $\psi_1^+$ as follows.
\begin{enumerate}[i.]
  \item $\psi_1^+(\lambda)=(\phi^+(\lambda))^+$ if $\lambda_\ell$ is odd and $\lambda_i$ are even for all $1\leq i\leq \ell-1$;
  \item $\psi_1^+(\lambda)=\phi^+_c(\phi^+(\lambda))$
      otherwise.
\end{enumerate}
Then we have both $\psi^+_1(\lambda)\in\mathcal{P}_\omega$ and $\psi^+_2(\lambda)\in\mathcal{P}_\omega$.
\end{lem}
\pf
For any given $\lambda=(\lambda_1,\ldots,\lambda_\ell)\in\mathcal{P}_{\omega}$,
by the definitions of $\mathcal{P}_\omega$ and $\phi^+_o$,
it is obvious  $\psi^+_2(\lambda)\in\mathcal{P}_\omega$.
For example, let $\lambda=(5,3,3,2)\in\mathcal{P}_\omega$,
then $\psi^+_2(\lambda)=(7,5,5,3)\in\mathcal{P}_\omega$.

Suppose that $\lambda_\ell$ is the only odd part in $\lambda$,
which means that all parts in $\phi^+(\lambda)$ are even.
Since  $(\phi^+(\lambda))^+=(\lambda_1,\lambda_2,\ldots,\lambda_{\ell-1},\lambda_{\ell}+1,0)$,
we  have $(\phi^+(\lambda))^+\in\mathcal{P}_\omega$.
As an example, let $\lambda=(6,6,4,3)\in\mathcal{P}_\omega$,
then $\psi^+_1(\lambda)=(6,6,4,4,0)\in\mathcal{P}_\omega$.

Otherwise,
it is easy to verify that there is at least one odd part in $\phi^+(\lambda)=(\lambda_1,\lambda_2,\ldots,\lambda_{\ell-1},\\\lambda_{\ell}+1)$.
If $\phi^+(\lambda)_\ell=\lambda_{\ell}+1$ is the only odd part,
then $\phi^+_c(\phi^+(\lambda))=(\lambda_1,\lambda_2,\ldots,\lambda_{\ell-1},\lambda_{\ell}/2+1,\lambda_{\ell}/2)$
so that $\phi^+_c(\phi^+(\lambda))\in\mathcal{P}_\omega$.
For example, let $\lambda=(8,8,4,2)\in\mathcal{P}_\omega$,
then $\phi^+(\lambda)=(8,8,4,3)$ and $\psi^+_1(\lambda)=(8,8,4,2,1)\in\mathcal{P}_\omega$.
Suppose the largest odd part in $\phi^+(\lambda)$ is $\lambda_i$ for some $1\leq i\leq \ell-1$.
Then we can assume that $\lambda_i=2b+1$
and $\lambda_\ell+1=a$ subject to $b\leq a-1$ since $\lambda\in\mathcal{P}_\omega$
and $\lambda_\ell=a-1$.
Thus $(\lambda_i-1)/2=b$
is the unique smallest part of $\phi^+_c(\phi^+(\lambda))$
and no odd part in $\phi^+_c(\phi^+(\lambda))$  exceeds $2b+1$,
implying that $\phi^+_c(\phi^+(\lambda))\in\mathcal{P}_\omega$.
For example, let $\lambda=(6,5,4,3)\in\mathcal{P}_\omega$,
then $\psi^+_1(\lambda)=(6,4,4,3,2)\in\mathcal{P}_\omega$.
\qed

\noindent{\it Combinatorial proof of \eqref{thmeq3}.}
We proceed to give a bijection between $\mathcal{P}_{\omega}(m,n)$ and $\mathcal{B}^1_{\omega}(m,n)$.
For a partition $\lambda=(\lambda_1,\ldots,\lambda_{m+1})\in\mathcal{P}_{\omega}(m,n)$,
let $\lambda^0=\lambda$, then for $i\geq1$,
we utilize the \emph{destructive} operator $\psi^-$ to construct $\lambda^{i}$ from $\lambda^{i-1}$ by
setting $\lambda^i=\psi^-(\lambda^{i-1})$ until $\ell(\lambda^i)=1$.
Denote by $t(\lambda)$ the last $i$ terminating the procedure
and   $d^\lambda=(d^\lambda_1,\ldots,d^\lambda_{t(\lambda)})$
the difference sequence recording the size difference  between $\lambda^{i-1}$ and $\lambda^{i}$,
that is, $d^\lambda_i=d_{\psi^-}(\lambda^{i-1})$ for  $1\leq i\leq t(\lambda)$.
By Lemma \ref{lem1},
we have $\lambda^i\in\mathcal{P}_\omega$ for $1\leq i\leq t(\lambda)$.
It can be easily checked that the number of 1's in $d^\lambda$ equals $m$
and $|\lambda^{t(\lambda)}|+\sum_{i=1}^{t(\lambda)}d^\lambda_i=n$.
As an example,
let $\lambda=(6,4,3,3,2)\in\mathcal{P}_\omega(4,18)$,
then the detailed steps of the procedure are listed in Table \ref{tabofdealg}.
\begin{table}
\begin{center}
\begin{tabular}{lc}
  \toprule
  partitions & size differences \\
  \midrule
  $\lambda^1=(6,5,4,2)$   & $d_{\psi^-}(\lambda^0)=1$ \\[3pt]
  $\lambda^2=(4,3,2,1)$   & $d_{\psi^-}(\lambda^1)=7$ \\[3pt]
  $\lambda^3=(4,3,2)$     & $d_{\psi^-}(\lambda^2)=1$ \\[3pt]
  $\lambda^4=(5,3)$       & $d_{\psi^-}(\lambda^3)=1$ \\[3pt]
  $\lambda^5=(3,2)$       & $d_{\psi^-}(\lambda^4)=3$ \\[3pt]
  $\lambda^6=(4)$         & $d_{\psi^-}(\lambda^5)=1$ \\[3pt]
  \bottomrule
\end{tabular}
\caption{detailed construction for  $\lambda=(6,4,3,3,2)$}\label{tabofdealg}
\end{center}
\end{table}
Hence we obtain $t(\lambda)=6$, $\lambda^6=(4)$ and $d^\lambda=(1,7,1,1,3,1)$.

\begin{claim}\label{claim1}
We have $d^\lambda_{t(\lambda)}=1$
and
for $1\leq i\leq t(\lambda)-1$,
if $d^\lambda_i=2k+1$ for some $k\geq1$,
then
$$\left|\left\{j\colon\,d^\lambda_{j}=1,\,i+1\leq j \leq t(\lambda)\right\}\right|=k.$$
\end{claim}

\noindent{\emph{Proof of Claim \ref{claim1}.}}
First noting that for $1\leq i\leq t(\lambda)$, $\ell(\lambda^{i-1})-\ell(\lambda^i)=1$
if and only if $d_{\psi^-}(\lambda^{i-1})=1$,
then by the constructing rules, it is evident that
$\ell(\lambda^{t(\lambda)})=1$ and $\ell(\lambda^{t(\lambda)-1})=2$,
which implies $d^\lambda_{t(\lambda)}=1$.
For $1\leq i\leq t(\lambda)-1$,
if $d_i^\lambda=d_{\psi^-}(\lambda^{i-1})=2k+1>1$,
we can deduce that $\psi^-(\lambda^{i-1})=\phi^-_o(\lambda^{i-1})$
and both $\lambda^{i-1}$ and $\lambda^{i}$ have $k+1$ parts by Lemma \ref{lem1}.
Hence we need to decrease the length of $\lambda^i$ from $k+1$
to 1 by iteratively acting $\psi^-$,
which means that  the number of $j$'s such that $d_j^\lambda$=1 for $i+1\leq j\leq t(\lambda)$ is exactly $k$.
\qed

\begin{figure}
  \centering
  \includegraphics[width=0.85\textwidth]{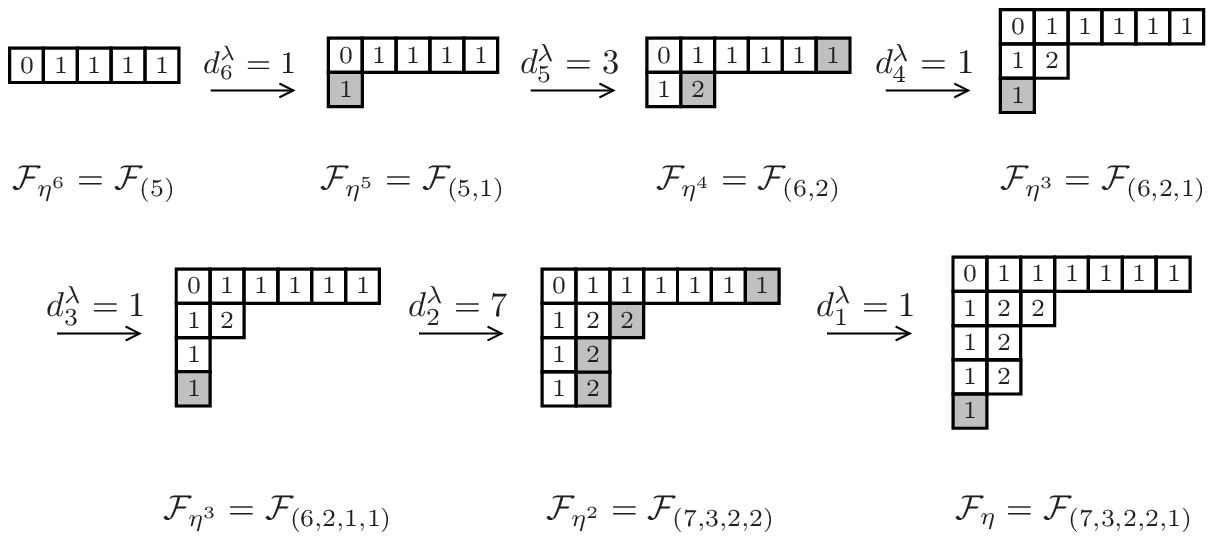}
  \caption{the procedure  of constructing $\f_\eta=\f_{(7,3,2,2,1)}$
  from $\lambda=(6,4,3,3,2)$,
  where the shaded boxes are those added at each step.}
  \label{opbijection}
\end{figure}

Now  we will construct the corresponding odd Ferrers graph $\f_\eta\in\mathcal{B}_\omega^1(m,n)$
from the partition $\lambda^{t(\lambda)}$
and the difference sequence $d^\lambda$.
Let $t=t(\lambda)$ and $\eta^{t}=(|\lambda^t|+1)$ be the partition with only one part $|\lambda^t|+1$,
then $\f_{\eta^{t}}\in\mathcal{B}_{\omega}$ is
of shape $\eta^{t}$ and size $|\lambda^t|$.
For $i$ ranges from $t-1$ to $0$,
we construct $\f_{\eta^{i}}$ from $\f_{\eta^{i+1}}$
and $d^\lambda_{i+1}$ as follows.
If $d^\lambda_{i+1}=1$,
we add a new row containing only one box below $\f_{\eta^{i+1}}$
and fill this box with 1,
which means $\ell(\f_{\eta^{i}})=\ell(\f_{\eta^{i+1}})+1$
and $|\f_{\eta^{i}}|=|\f_{\eta^{i+1}}|+1=|\lambda^{i}|$.
If $d^\lambda_{i+1}=2k+1$ for some $k\geq1$,
we  add one box at the end of each row of $\f_\eta^{i+1}$,
and fill the added box in the first row with 1
and the rest of the boxes with 2's.
It can be seen that $\ell(\f_{\eta^{i}})=\ell(\f_{\eta^{i+1}})$
and by  Claim \ref{claim1}, $|\f_{\eta^{i}}|=|\f_{\eta^{i+1}}|+(2k+1)=|\lambda^{i}|$.
Finally, setting $\f_\eta=\f_{\eta^0}$,
since there are exactly $m$ 1's in the difference sequence $d^\lambda$,
we have $\f_\eta\in\mathcal{B}_\omega^1(m,n)$ as desired.
As the above example, $t=6$, $\lambda^6=4$ and $d^\lambda=(1,7,1,1,3,1)$,
the corresponding $\f_\eta$ is constructed in Figure \ref{opbijection}.

To complete the proof, it remains to construct the corresponding partition $\lambda\in\mathcal{P}_{\omega}(m,n)$
from a given odd Ferrers graph $\f_\eta\in\mathcal{B}^1_{\omega}(m,n)$.
Let $\f_{\eta^0}=\f_\eta$,
then for $i\geq 1$,
we will construct $\f_{\eta^{i}}$ from  $\f_{\eta^{i-1}}=\f_{(\eta^{i-1}_1,\ldots,\eta^{i-1}_{\ell})}$
and difference sequence $h^{\f_\eta}=(h^{\f_\eta}_1,h^{\f_\eta}_2,\ldots)$
as follows.

\noindent {Case 1}:
Let $\eta^{i}=\phi^{*}(\eta^{i-1})$
if $\eta^{i-1}_{\ell}\geq2$,
then the corresponding odd Ferrers graph $\f_{\eta^i}$ is obtained by
deleting the last box of each row  of $\f_{\eta^{i-1}}$.
Let $h^{\f_\eta}_i$ be the sum of the deleted numbers.

\noindent {Case 2}:
Let $\eta^{i}=\phi^{-}(\eta^{i-1})$ if $\eta^{i-1}_\ell=1$,
then the corresponding odd Ferrers graph $\f_{\eta^i}$ is obtained
by deleting the single box with a 1 in the last row  of $\f_{\eta^{i-1}}$.
As in Case 1, let $h^{\f_\eta}_i$ be the sum of the deleted numbers, which equals 1.

\noindent This process will not stop until $\ell(\f_{\eta^i})=1$ for some $i\geq1$.
Denoting this $i$ by $r(\f_\eta)$,
we can deduce that  $h^{\f_\eta}_i$ is odd
for all $1\leq i\leq r(\f_\eta)$.
Moreover, the number of $1$'s in the difference sequence $h^{\f_\eta}$ is $m$ and
$|\f_{\eta^{r}}|+\sum_{i=1}^r h^{\f_\eta}_i=n$, where $r=r(\f_\eta)$.
Following a similar analysis to that in Claim \ref{claim1},
for the difference sequence $h^{\f_\eta}$,
we also have the following.

\begin{claim}\label{claim2}
We have $h^{\f_\eta}_{r}=1$
and for $1\leq i\leq r-1$,
if $h^{\f_\eta}_i=2k+1$ for some $k\geq1$,
then
$$\left|\left\{j\colon\,h^{\f_\eta}_{j}=1,\,i+1\leq j \leq r\right\}\right|=k.$$
\end{claim}

According to the difference sequence $h^{\f_\eta}$ and $\f_{\eta^r}$,
we can use \emph{constructive} operators $\psi^+_1$ and $\psi^+_2$ to construct $\lambda$.
Firstly, let $\lambda^{r}=(|\eta^{r}|-1)$ be the partition with only one part $|\eta^{r}|-1$,
so that  $\lambda^r\in\mathcal{P}_{\omega}$
and $|\lambda_r|=|\f_{\eta^r}|.$
For $r-1\geq i\geq 0$,
suppose that $\lambda^{i+1}$ has already been generated,
we generate $\lambda^i$ by letting
$\lambda^{i}=\psi^+_1(\lambda)$
if $h^{\f_\eta}_{i+1}=1$
or $\lambda^{i}=\psi^+_2(\lambda)$
if $h^{\f_\eta}_{i+1}=2k+1$ for some $k\geq1$.
Since the length of $\lambda^{i}$  increases 1
if and only if $h^{\f_\eta}_{i+1}=1$,
then by Claim \ref{claim2},
we obtain $|\lambda^i|=|\lambda^{i+1}|+h^{\f_\eta}_{i+1}$.
Therefore, by setting $\lambda=\lambda^0$ and Lemma \ref{lem2},
we see that $\lambda\in\mathcal{P}_\omega(m,n)$.
\qed

\begin{example}
Based on the above bijection,
we give the one-to-one correspondence between all $\lambda\in\mathcal{P}_\omega(3,15)$
and $\f_\eta\in\mathcal{B}^1_\omega(3,15)$ in Table \ref{tabcorop}.
\begin{table}
\begin{center}
\begin{tabular}{cc|cc}
  \toprule
  $\lambda\in\mathcal{P}_\omega(3,15)$ & $\f_\eta\in\mathcal{B}^1_\omega(3,15)$ & $\lambda\in\mathcal{P}_\omega(3,15)$ & $\f_\eta\in\mathcal{B}^1_\omega(3,15)$ \\
  \midrule
  $(12,2,1,0)$ & $\f_{(11,2,1,1)}$ & $(8,3,3,1)$ & $\f_{(7,2,2,2)}$ \\[3pt]
  $(10,4,1,0)$ & $\f_{(7,4,1,1)}$ & $(6,4,4,1)$ & $\f_{(5,3,2,2)}$\\[3pt]
  $(8,6,1,0)$  & $\f_{(9,3,1,1)}$ & $(6,4,3,2)$ & $\f_{(7,3,2,1)}$\\[3pt]
  $(10,2,2,1)$ & $\f_{(9,2,2,1)}$ & $(5,5,3,2)$ & $\f_{(5,5,1,1)}$ \\[3pt]
  $(8,4,2,1)$  & $\f_{(5,4,2,1)}$ & $(5,4,4,2)$ & $\f_{(3,3,3,2)}$\\[3pt]
  $(6,6,2,1)$  & $\f_{(5,3,3,1)}$ & $(4,4,4,3)$ & $\f_{(13,1,1,1)}$\\[3pt]
  \bottomrule
\end{tabular}
\end{center}
\caption{correspondence between $\mathcal{P}_\omega(3,15)$ and $\mathcal{B}^1_\omega(3,15)$}
\label{tabcorop}
\end{table}
\end{example}

\subsection{Bijection between $\mathcal{P}_{\nu}(m,n)$ and $\mathcal{B}_\nu^1(m,n)$}
Recall that $\mathcal{P}_{\nu}$ is the
set of partitions with distinct parts
which may includes 0, satisfying that
the odd parts are less than twice  the smallest part.
Before giving the bijective proof of \eqref{eqthm2}, we need the following two lemmas.

\begin{lem}\label{lem3}
Given any  partition $\lambda=(\lambda_1,\ldots,\lambda_\ell)\in\mathcal{P}_{\nu}$,
define the \emph{destructive} operator $\rho^-$ as follows.
\begin{enumerate}[i.]
  \item $\rho^-(\lambda)=\phi^-_e(\lambda^-)$ if $\lambda_{\ell}=0$;
  \item $\rho^-(\lambda)=\phi^-_e(\phi^-_c(\lambda))$ if $\lambda_{\ell}\geq 1$ and $\lambda_{\ell-1}=\lambda_{\ell}+1$;
  \item $\rho^-(\lambda)=\phi^-_o(\lambda)$ if $\lambda_{\ell}\geq1$ and $\lambda_{\ell-1}\geq \lambda_{\ell}+2$.
\end{enumerate}
Then we have $\rho^-(\lambda)\in\mathcal{P}_\nu$.
\end{lem}
\pf
According  to the quantitative relation between $\lambda_\ell$ and $\lambda_{\ell-1}$,
and the actions of $\phi^-_o$ and $\phi^-_e$,
it can be directly checked that
all parts in $\rho^-(\lambda)$ are distinct. Hence, we only need to examine the restriction on odd parts.

If $\lambda_\ell=0$,
by the definition of $\mathcal{P}_\nu$,
all parts of $\lambda$ are even and $\lambda_{\ell-1}\geq2$.
Since $\phi^-_e(\lambda^-)=(\lambda_1-2,\lambda_2-2,\ldots,\lambda_{\ell-1}-2)$,
it follows that all parts  of $\phi^-_e(\lambda^-)$ are even
and $\lambda_{\ell-1}-2\geq0$.
Hence we have $\phi^-_e(\lambda^-)\in\mathcal{P}_\nu$.
For example, let $\lambda=(12,8,4,2,0)\in\mathcal{P}_\nu$, then $\rho^-(\lambda)=(10,6,2,0)\in\mathcal{P}_\nu$.

If $\lambda_\ell\geq1$
and $\lambda_{\ell-1}=\lambda_{\ell}+1$,
assuming  $\lambda_\ell=a$ and $\lambda_0=\infty$,
it is clear that there exists a unique $i_0$,
with $1\leq i_0\leq \ell-1$, such that $\lambda_{i_0}<2a+1<\lambda_{i_0-1}$.
Thus $\phi^-_c(\lambda)=(\lambda_1,\lambda_2,\ldots,\lambda_{i_0-1},2a+1,\lambda_{i_0},\ldots,\lambda_{\ell-2})$.
Since $\lambda\in\mathcal{P}_\nu$,
we know that all odd parts in $\lambda$ are less than $2a$,
which implies that the largest odd part in $\phi^-_c(\lambda)$ is $2a+1$.
By the definition of $\phi^-_e$, we obtain
$\phi^-_e(\phi^-_c(\lambda))=(\lambda_1-2,\lambda_2-2,\ldots,\lambda_{i_0-1}-2,2a-1,\lambda_{i_0}-2,\ldots,\lambda_{\ell-2}-2)$,
which leads to the largest odd part in $\phi^-_e(\phi^-_c(\lambda))$ being $2a-1$.
If $i_0=\ell-1$ then we are done,
otherwise the smallest part in $\phi^-_e(\phi^-_c(\lambda))$ is $\lambda_{\ell-2}-2$.
Note that the partition $\lambda$ is distinct and $\lambda_\ell=a$,
so it is clear that $\lambda_{\ell-2}-2\geq\lambda_\ell+2-2\geq a$.
Hence we deduce that $\phi^-_e(\phi^-_c(\lambda))\in\mathcal{P}_\nu$.
For example, let $\lambda=(10,8,7,5,4)\in\mathcal{P}_\nu$,
then $\rho^-(\lambda)=(8,7,6,5)\in\mathcal{P}_\nu$.

If $\lambda_\ell\geq1$
and $\lambda_{\ell-1}\geq\lambda_{\ell}+2$,
we have $\phi^-_o(\lambda)=(\lambda_1-2,\lambda_2-2,\ldots,\lambda_{\ell-1}-2,\lambda_{\ell}-1)$.
For $1\leq i\leq \ell-1$, suppose that $\lambda_i$ is odd,
then by the definition of $\mathcal{P}_{\nu}$,
we conclude that $\lambda_i\leq 2\lambda_\ell-1$,
which implies that
$\lambda_i-2\leq 2\lambda_\ell-3=2(\lambda_\ell-1)-1$, thus $\phi^-_o(\lambda)\in\mathcal{P}_{\nu}$.
For example, let $\lambda=(8,6,5,3)\in\mathcal{P}_\nu$,
then $\rho^-(\lambda)=(6,4,3,2)\in\mathcal{P}_\nu$.
\qed

\begin{lem}\label{lem4}
Given any  partition $\lambda=(\lambda_1,\ldots,\lambda_\ell)\in\mathcal{P}_{\nu}$,
define the \emph{constructive} operators $\rho^+_2$ by $\rho^+_2(\lambda)=\phi^+_o(\lambda)$ and $\rho^+_1$ as follows.
\begin{enumerate}[i.]
  \item $\rho^+_1(\lambda)=(\phi^+_e(\lambda))^+$
  if $\lambda_{i}$ is even for all $1\leq i\leq \ell$;
  \item $\rho^+_1(\lambda)=\phi^+_c(\phi^+_e(\lambda))$
  if $\lambda_{i}$ is odd for some $1\leq i\leq \ell$.
\end{enumerate}
Then we have both $\rho^+_1(\lambda)\in\mathcal{P}_\nu$ and $\rho^+_2(\lambda)\in\mathcal{P}_\nu$.
\end{lem}

The proof of Lemma \ref{lem4} is routine and similar to Lemma \ref{lem1}, Lemma \ref{lem2} and Lemma \ref{lem3}
so we omit the proof here and just present some examples.
Let $\lambda=(10,8,7,6,5)\in\mathcal{P}_\nu$,
then $\rho^+_2(\lambda)=(12,10,9,8,6)\in\mathcal{P}_\nu$.
Let $\lambda=(10,8,4,2)\in\mathcal{P}_\nu$,
then $\rho^+_1(\lambda)=(\phi^+_e(\lambda))^+=(12,10,6,4,0)\in\mathcal{P}_\nu$
and let $\lambda=(10,8,7,4)\in\mathcal{P}_\nu$,
then $\rho^+_1(\lambda)=\phi^+_c(\phi^+_e(\lambda))=(12,10,6,5,4)\in\mathcal{P}_\nu$.

\noindent{\it Combinatorial proof of \eqref{eqthm2}.}
By \eqref{pnu} and \eqref{comtooprow},
we will construct the odd Ferrers graph $\f_\eta=\f_{(\eta_1,\ldots,\eta_{m+1})}\in\mathcal{B}^1_{\nu}(m,n)$
from a partition $\lambda=(\lambda_1,\ldots,\lambda_{m+1})\in\mathcal{P}_{\nu}(m,n)$
by the following procedure.
To this end, set $\lambda^0=\lambda$
and suppose that for $i\geq 1$ the partition $\lambda^{i-1}$ has been constructed.
Then we will continue constructing the partition $\lambda^i$ by
setting $\lambda^i=\rho^-(\lambda^{i-1})$
until $\ell(\lambda^i)=1$.
Denote by $t(\lambda)$ the $i$ satisfying $\ell(\lambda^i)=1$ and
$d^\lambda=(d^\lambda_1,\ldots,d^\lambda_{t(\lambda)})$
the difference sequence with $d^\lambda_i=d_{\rho^-}(\lambda^{i-1})$.
Thus we have $\lambda_i\in\mathcal{P}_\nu$ for each $1\leq i\leq t(\lambda)$ by Lemma \ref{lem3}.
It is easy to see that
the number of even numbers in $d^\lambda$ is $m$
and $|\lambda^t|+\sum_{i=1}^{t(\lambda)}d^\lambda_i=n$.
For example, let $\lambda=(10,8,5,4,3)\in\mathcal{P}_{\nu}(4,30)$,
from Table \ref{tabofdpdealg}
we have $t(\lambda)=5$, $\lambda^5=(3)$ and $d^\lambda=(8,7,6,4,2)$.
\begin{table}
\begin{center}
\begin{tabular}{lc}
  \toprule
  partitions & size differences \\
  \midrule
  $\lambda^1:(8,6,5,3)$   & $d_{\rho^-}(\lambda^0)=8$ \\[3pt]
  $\lambda^2:(6,4,3,2)$   & $d_{\rho^-}(\lambda^1)=7$ \\[3pt]
  $\lambda^3:(4,3,2)$     & $d_{\rho^-}(\lambda^2)=6$ \\[3pt]
  $\lambda^4:(3,2)$       & $d_{\rho^-}(\lambda^3)=4$ \\[3pt]
  $\lambda^5:(3)$       & $d_{\rho^-}(\lambda^4)=2$ \\[3pt]
  \bottomrule
\end{tabular}
\caption{detailed construction for $\lambda=(10,8,5,4,3)$}\label{tabofdpdealg}
\end{center}
\end{table}

\begin{claim}\label{claim3}
We have $d^\lambda_{t(\lambda)}=2$
and
for $1\leq i\leq t(\lambda)-1$,
if $d^\lambda_i=2k+1$ for some $k\geq1$,
then
$$\left|\left\{j\colon\,d^\lambda_{j}\mbox{ is even},\,i+1\leq j \leq t(\lambda)\right\}\right|=k.$$
\end{claim}

\noindent{\emph{Proof of Claim \ref{claim3}}.}
By the definition of the operator $\rho^-$,
we know that for $1\leq i\leq t(\lambda)$,
$\ell(\lambda^{i-1})=\ell(\lambda^i)+1=k+1$
if and only if
$d^\lambda_i=2k$ for some $k\geq1$.
Thus $d^{\lambda}_{t(\lambda)}=2$
since $\ell(\lambda^{t(\lambda)})=1$
and $\ell(\lambda^{t(\lambda)-1})=2$
by the terminating condition.
For $1\leq i\leq t(\lambda)-1$,
if $d_i^\lambda=2k+1>1$,
we can easily deduce that both $\lambda^{i-1}$ and $\lambda^{i}$ have $k+1$ parts by Lemma \ref{lem3}.
Since the length of $\lambda^i$ is decreased
from $k+1$ to 1,
there are exactly $k$ even $d_j^\lambda$'s
in the set $\{d_{i+1}^\lambda,\ldots,d_{t(\lambda)}^\lambda\}$.
\qed

\begin{figure}
  \centering
  \includegraphics[width=0.85\textwidth]{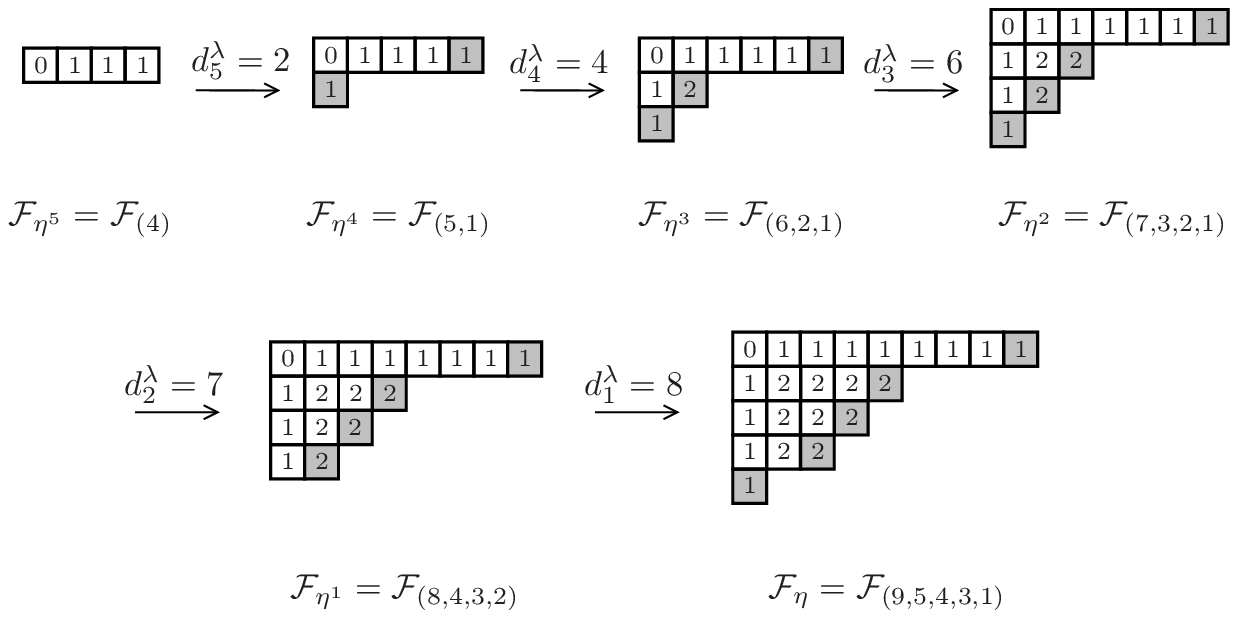}
  \caption{the procedure of constructing $\f_\eta=\f_{(9,5,4,3,1)}$
  from $\lambda=(10,8,5,4,3)$,
  where the shaded boxes are those added  at each step.}
  \label{dpbijection}
\end{figure}

Setting $t=t(\lambda)$, we will construct $\f_\eta$ by $\lambda^t$ and $d^\lambda$ as follows.
Let $\eta^t$ be the partition with only one part $|\lambda^t|+1$,
then we have
$\f_{\eta^t}\in\mathcal{B}_\nu$
and $|\f_{\eta^t}|=|\lambda^t|$.
For $t-1 \geq i\geq 0$,
suppose that we have obtained $\f_{\eta^{i+1}}$.
If $d^\lambda_{i+1}$ is even,
we construct $\f_{\eta^i}$ by first adding a new box at the end of each row of  $\f_{\eta^{i+1}}$
and a new row with only one box
under the bottom of $\f_{\eta^{i+1}}$,
then filling the added boxes in the first and new rows by 1's and the other boxes by 2's.
If $d^\lambda_{i+1}$ is odd,
we add one box at the end of each row of $\f_{\eta^{i+1}}$,
and fill the added boxes with 2's but the box in the first row with 1.
Note that these manipulations ensure that $\f_{\eta^i}\in\mathcal{B}_\nu$
and $\ell(\f_{\eta^i})=\ell(\f_{\eta^{i+1}})+1$ if $d^\lambda_{i+1}$ is even,
$\ell(\f_{\eta^i})=\ell(\f_{\eta^{i+1}})$ if $d^\lambda_{i+1}$ is odd.
By Claim \ref{claim3},
it can be seen $\ell(\f_{\eta^{i+1}})=\left\lceil(d^\lambda_{i+1}+1)/2\right\rceil$
and $|\f_{\eta^{i}}|=|\f_{\eta^{i+1}}|+|d^\lambda_{i+1}|$
for $0\leq i\leq t-1$.
Thus by letting $\f_\eta=\f_{\eta^0}$,
we arrive at
$\f_\eta\in\mathcal{B}^1_{\nu}(m,n)$.
As the above example
for $\lambda=(10,8,5,4,3)$,
the corresponding $\f_\eta=\f_{(9,5,4,3,1)}\in\mathcal{P}_\nu(4,30)$ is generated in Figure \ref{dpbijection}.

To complete the proof,
we need to construct the corresponding partition $\lambda\in\mathcal{P}_{\nu}(m,n)$
from a given odd Ferrers graph $\f_\eta\in\mathcal{B}^1_{\nu}(m,n)$.
Let $\f_{\eta^0}=\f_\eta$.
For $i\geq 1$,
we construct $\f_{\eta^i}$  by letting
$\f_{\eta^i}=\f_{(\phi^{*}(\eta^{i-1}))}$
until $\ell(\f_{\eta^i})=1$ for some $i$.
Denote by $r(\f_\eta)$ the $i$ such that  $\ell(\f_{\eta^i})=1$
and $h^{\f_\eta}=(h_1^{\f_\eta},\ldots,h_{r(\f_\eta)}^{\f_\eta})$
the difference sequence,
where $h_i^{\f_\eta}=|\f_{\eta^{i-1}}|-|\f_{\eta^i}|$
for $1\leq i\leq r(\f_\eta)$.
Now we can generate $\lambda$
by using \emph{constructive} operators $\rho^+_1$ and $\rho^+_2$ as follows.
Let $r=r(\f_\eta)$ and
$\lambda^r=(|\f_{\eta^r}|-1)$ be the partition with only one part $|\f_{\eta^r}|-1$,
then it follows that $\lambda^r\in\mathcal{P}_\nu$ and
$|\lambda^r|=|\f_{\eta^r}|$.
For $i$ ranging from $r-1$ to $0$,
let $\lambda^i=\rho^+_1(\lambda)$
if $h^{\f_\eta}_{i+1}$ is even
and let $\lambda^i=\rho^+_2(\lambda)$
if $h^{\f_\eta}_{i+1}$ is odd.
Therefore, by Lemma \ref{lem4},
we obtain that
$\lambda=\lambda^0\in\mathcal{P}_\nu(m,n)$.
\qed

\begin{example}
To conclude this section, we present the one-to-one correspondence between all $\lambda\in\mathcal{P}_\nu(4,30)$ and $\f_\eta\in\mathcal{B}_\nu^1(4,30)$ in Table \ref{tabcordp}.
\begin{table}
\begin{center}
\begin{tabular}{cc|cc}
  \toprule
  $\lambda\in\mathcal{P}_\nu(4,30)$ & $\f_\eta\in\mathcal{B}^1_\nu(4,30)$ & $\lambda\in\mathcal{P}_\nu(4,30)$ & $\f_\eta\in\mathcal{B}^1_\nu(4,30)$ \\
  \midrule
  $(18,6,4,2,0)$ & $\f_{(15,4,3,2,1)}$ & $(12,8,6,4,0)$ & $\f_{(7,6,4,3,1)}$ \\[3pt]
  $(16,8,4,2,0)$ & $\f_{(11,6,3,2,1)}$ & $(10,8,6,4,2)$ & $\f_{(7,5,4,3,2)}$\\[3pt]
  $(14,10,4,2,0)$  & $\f_{(9,7,3,2,1)}$ & $(12,6,5,4,3)$ & $\f_{(11,5,4,2,1)}$\\[3pt]
  $(14,8,6,2,0)$ & $\f_{(9,6,4,2,1)}$ & $(10,8,5,4,3)$ & $\f_{(9,5,4,3,1)}$ \\[3pt]
  $(12,10,6,2,0)$  & $\f_{(7,6,5,2,1)}$ & $(8,7,6,5,4)$ & $\f_{(13,5,3,2,1)}$\\[3pt]
  \bottomrule
\end{tabular}
\caption{correspondence between $\mathcal{P}_\nu(4,30)$ and $\mathcal{B}^1_\nu(4,30)$}
\label{tabcordp}
\end{center}
\end{table}
\end{example}

\section{Further remarks}\label{secofrem}
Although the trivariate generalizations $\omega(y,z;q)$ and $\nu(y,z;q)$ of the mock theta functions $\omega(q)$ and $\nu(q)$
have been constructed and studied,
the identities related to $\omega(y,z;q)$ and $\nu(y,z;q)$
that are analogous to Theorem \ref{thm1} still remain mysterious.
For the sake of analytic approach,
one may utilize the arithmetic properties of $\bar{\nu}(\alpha,z;q)$ and $\bar{\omega}(\alpha,z;q)$
investigated in \cite{choi},
since $\nu(y,z;q)$ and $\omega(y,z;q)$ are related to $\bar{\nu}(\alpha,z;q)$ and $\bar{\omega}(\alpha,z;q)$ by \eqref{nutobarnu} and \eqref{omegatobaromega}, respective.
From the perspective of combinatorics,
the variable $y$ in $\omega(y,z;q)$ or $\nu(y,z;q)$
is related to the columns of the corresponding odd Ferrers graphs,
but from the constructive algorithms  given in Section \ref{bjisect}, it is difficult to determine which partition statistic in $\mathcal{P}_\omega$
or $\mathcal{P}_\nu$ can be reflected by $y$.
Thus the simulation of Theorem \ref{thm1} on generalized trivariate mock theta functions $\omega(y,z;q)$ and $\nu(y,z;q)$
is desired by either combinatorial or analytic methods.

\vspace{0.5cm}

\noindent{\bf Acknowledgements.}
The authors appreciate the referee for his/her helpful comments which improve the quality of this manuscript.
This work is supported by  Doctor Scientific Research Foundation of Chongqing University of Posts and Telecommunications (Grant No. A2017-123).

\end{document}